\documentclass[10pt]{article}
\usepackage[english]{babel}
\usepackage{natbib}
\usepackage{url}
\usepackage{inputenc}
\usepackage{amsfonts}
\usepackage{amsmath}
\usepackage{empheq}
\usepackage{graphicx}
\graphicspath{{figures/}}
\usepackage{fancyhdr}
\usepackage{bbold}
\usepackage{amssymb}
\usepackage{vmargin}
\usepackage{array}
\usepackage{amsthm}
\usepackage[ruled,vlined]{algorithm2e}
\usepackage{hyperref}
\usepackage{caption}
\usepackage{float}
\usepackage{fourier-orns}
\usepackage{listings}

\newtheorem{theorem}{Theorem}[section]
\newtheorem{corollary}[theorem]{Corollary}
\newtheorem{lemma}[theorem]{Lemma}
\newtheorem{remark}[theorem]{Remark}
\theoremstyle{definition}

\newtheorem{prop}[theorem]{Proposition}

\counterwithin{equation}{section}

\usepackage[shortlabels]{enumitem}
\mathtoolsset{showonlyrefs=true}
\usepackage[scr=boondoxo,scrscaled=1.05]{mathalfa}

\newcommand{\om}{{\Omega}}

\newcommand{\bydef}{\stackrel{\mbox{\tiny\textnormal{\raisebox{0ex}[0ex][0ex]{def}}}}{=}}

\usepackage[dvipsnames]{xcolor}

\title{Stability of patterns in variational systems \\ with continuous symmetries: \\ a computer-assisted approach}

\author{
Jan Bouwe van den Berg
\footnote{VU Amsterdam, Department of Mathematics, De Boelelaan 1081, 1081 HV Amsterdam, The Netherlands. {\tt janbouwe@few.vu.nl}}
\and
Dominic Blanco
\footnote{Rutgers University, Department of Mathematics, 110 Frelinghuysen Road, Piscataway, NJ, 08854. {\tt dominic.blanco@rutgers.edu}}
}
\begin{document}
\maketitle
\begin{abstract}
In variational models of pattern formation, periodic solutions are never isolated: since the energy does not change under translations, every shifted copy of a pattern is again a solution. This is an obstacle for computer-assisted proofs, which need isolated solutions. Imposing discrete symmetries removes the obstacle, but restricts the patterns that can be studied and only gives stability against symmetric perturbations. We develop a computer-assisted method that deals with the translations directly. We fix the position of the pattern with phase conditions and add unfolding parameters. They act as Lagrange multipliers, and translation invariance of the energy ensures that they vanish at a solution, under a condition that is verified rigorously by the computer. We then derive the stability of the pattern from the computer-assisted proof. We obtain the number of unstable directions and show that the only neutral directions are the translations. We apply the method to the Swift-Hohenberg equation in two dimensions and prove the existence of stable and unstable patterns without imposed symmetry. Since the linear algebra does not depend on the number of translation directions, the approach extends to higher dimensions.

\end{abstract}

\section{Introduction}

Periodic patterns such as stripes and hexagonal arrays, but also less regular ones such as labyrinths and disordered arrangements of spots, appear throughout physics, chemistry and biology. In many models they are critical points of an energy, and the stable ones are local minimizers. Numerical simulations readily produce such patterns, but proving that a computed pattern corresponds to a true solution, and deciding whether it is stable, requires more. Computer-assisted proofs can do both, but translation invariance stands in their way, since a proof of this kind singles out one solution and every translate of a pattern is again a solution. The usual remedy is to impose a discrete symmetry, which restricts the patterns that can be studied and the perturbations against which stability is shown. In this paper we develop a computer-assisted method that handles the translations directly.

We consider stationary solutions of gradient flows whose energy is invariant under translations. To demonstrate our approach, we apply it to the Swift-Hohenberg (SH) equation~\cite{swift_hohenberg_1977}, a canonical model for pattern formation~\cite{cross_hohenberg_1993},
\begin{equation}\label{eq:swift-hohenberg-2d}
    \partial_t u = -(I+\Delta)^2 u - \mu u + \gamma u^2 - u^3,
\end{equation}
in two spatial dimensions on a \emph{periodic} domain $\om \bydef [-d_1,d_1] \times [-d_2,d_2]$, where $\Delta \bydef \partial_1^2 + \partial_2^2$ is the Laplacian. The SH equation is the gradient flow of the energy
\begin{equation}\label{eq:energy}
  \mathcal{E}(u) \bydef \frac{1}{|\om|} \int_{\om} \frac{1}{2} \bigl((I+\Delta)u\bigr)^2 + \frac{\mu}{2} u^2 - \frac{\gamma}{3} u^3 + \frac{1}{4} u^4 \, dx .
\end{equation}
That is, it reads $\partial_t u = -F(u)$, where
\[
    F(u) \bydef (I + \Delta)^2u + \mu u -\gamma u^2 + u^3
\]
is the $L^2$-gradient of $\mathcal{E}$, see Section~\ref{sec:variational}. The stationary solutions are the zeros of $F$.

Periodic solutions of the SH equation have been established with computer-assisted proofs~\cite{lessard_2,lessard_1}, for instance with hexagonal and triangular symmetry in~\cite{dominic_sh_periodic}. Localized patterns have been studied analytically and numerically~\cite{ladder2009Beck,jason_spot_paper,jason_ring_paper,squareSakaguchi_1997}; see~\cite{jason_review_paper} for a review. In the plane they come in families of translates and rotations, and the computer-assisted proofs in~\cite{symmetry_blanco_cadiot} isolate them by imposing dihedral symmetry.

\begin{figure}[t]
\centering\includegraphics[width=0.65\linewidth]{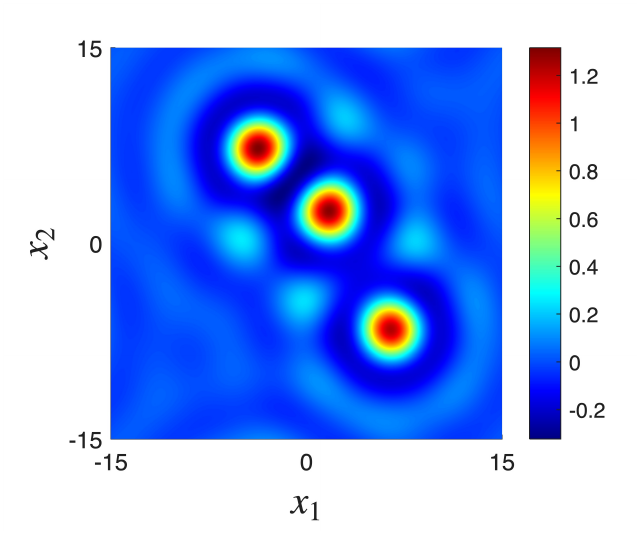}
\caption{The approximate solution $u_0$ of Theorem~\ref{th:intro}.}
 \label{fig:intro}
 \end{figure}

In this paper we consider periodic patterns without any imposed symmetry. We work on a periodic domain where the only continuous symmetries are the translations. The linearization $DF(\tilde{u})$ at a stationary solution $\tilde{u}$ has discrete spectrum, so that its number of negative eigenvalues and its nullspace are well defined. Stability is then with respect to perturbations with the same periods. We impose no symmetry because patterns found in simulations often have none, and because only then is stability shown against all such perturbations. A typical result is the following.
\begin{theorem}\label{th:intro}
Let $\mu = 0.25$, $\gamma = 1.6$ and $d_1 = d_2 = 15$, and let $u_0$ be the approximate solution described by finitely many Fourier coefficients (see Section~\ref{sec:results}) and shown in Figure~\ref{fig:intro}. There exists a stationary solution $\tilde{u}$ of~\eqref{eq:swift-hohenberg-2d} with $\|\tilde{u} - u_0\|_{\ell^1_{\nu}} \leq 4 \times 10^{-9}$, and $\tilde{u}$ is locally unique up to translations. The eigenvalue $0$ of $DF(\tilde{u})$ has multiplicity two, with eigenfunctions $\partial_1 \tilde{u}$ and $\partial_2 \tilde{u}$ coming from the translation invariance, and $DF(\tilde{u})$ has one negative eigenvalue. In particular, $\tilde{u}$ is unstable.
\end{theorem}

Our proofs are computer-assisted proofs (CAPs), which combine analysis with validated numerics in interval arithmetic; see~\cite{nakao_plum_watanabe_2019,gomez_serrano_2019} for overviews in the PDE context. They rely on a Newton-Kantorovich type contraction argument in the form of radii polynomials~\cite{day_lessard_mischaikow_2007} (see Theorem~\ref{th:radii-polynomial}). We compute an approximate solution $u_0$ and an approximate inverse $A$ of $DF(u_0)$, and show that $T(u) \bydef u - AF(u)$ is a contraction on a ball around $u_0$. Its unique fixed point $\tilde{u}$ in this ball is a zero of $F$. The stability of $\tilde{u}$ is determined by the eigenvalues of $DF(\tilde{u})$, which can be enclosed rigorously with Gershgorin type arguments and interval arithmetic (see e.g.~\cite{orrsommerfeld_eq,Jean-PhilippeLessard2025DiscreteandContinuousDynamicalSystems-B,rump2006}).

Translation invariance prevents a direct application of this argument. If $\tilde{u}$ is a stationary solution, then so is $\tilde{u}(\cdot + s)$ for every $s \in \mathbb{R}^2$. Solutions are thus not isolated and a naive contraction argument cannot succeed. As mentioned above, imposing a discrete symmetry removes this problem, as in the proofs for the SH equation. Other examples are the computer-assisted proofs for the Kuramoto-Sivashinsky equation in~\cite{zgliczynski_mischaikow_2001,wilczak_zgliczynski_2020}, which restrict to odd functions, and~\cite{JB_symmetries_1,JB_symmetries_2}, where space group symmetries are enforced in Fourier series and applied to the Ohta-Kawasaki problem in three dimensions. The case in which a continuous symmetry remains was left for future work in~\cite{JB_symmetries_2}. We treat it here, using full Fourier series.

The main contribution of this paper is a computer-assisted method for periodic patterns in variational systems that proves both existence and stability without imposed symmetry. Adding phase conditions and unfolding parameters has been used as a strategy previously (see e.g.~\cite{symmetry_blanco_cadiot,jay_unfold1,jay_unfold4,jay_unfold3,jay_unfold2}). The novelty here is that our formulation respects the variational structure, so that stability information is retained.
The method is simple: compared with a standard computer-assisted proof, it adds two phase conditions with their unfolding parameters and the invertibility check of one $2 \times 2$ matrix. The number of unstable directions then follows from the bounds of the existence proof, and one additional estimate gives the spectral gap.

We isolate the solution with two phase conditions, which fix the position of the solution among its translates, and two unfolding parameters $\beta_1,\beta_2$, which keep the number of equations and unknowns equal. The augmented system is again variational, with $\beta_1,\beta_2$ as Lagrange multipliers for the phase conditions (Section~\ref{sec:fourier-setting}). A zero of the augmented system gives a zero of the original $F$ of interest only if $\beta_1 = \beta_2 = 0$. We show that translation invariance of the energy indeed forces $\beta_1$ and $\beta_2$ to vanish, provided a $2 \times 2$ matrix is invertible, which we check with interval arithmetic (Section~\ref{sec:beta-zero}).

For the stability we need the spectrum of $DF(\tilde{u})$, whereas the contraction argument controls the derivative of the augmented system. This derivative is $DF$ bordered by two extra rows and columns, which come from the phase conditions and the unfolding parameters. We crucially choose these to contain $\partial_1 u_0$ and $\partial_2 u_0$, corresponding to the infinitesimal translations of the approximate solution. If we replace these by the derivatives $\partial_1 \tilde{u}$ and $\partial_2 \tilde{u}$ of the \emph{exact} solution, the linear algebra of bordered operators (Section~\ref{sec:stability-2D}) shows that the resulting operator has the same eigenvalues as $DF(\tilde{u})$, except that the double eigenvalue $0$ is replaced by two positive and two negative eigenvalues. 
Since $\tilde{u}$ is known only up to the error bound, we use one additional estimate, which shows that the approximate inverse from the existence proof is also an approximate inverse of the operator with the exact border
(Theorem~\ref{th:radii-polynomial-q} and Section~\ref{sec:Zq}). From this we obtain the number of negative eigenvalues of $DF(\tilde{u})$ by counting the negative eigenvalues of a \emph{finite} matrix, and we show that $0$ is an eigenvalue of multiplicity exactly two (Lemma~\ref{lem:eigenvalues-of-Q-2d} and Corollary~\ref{lem:homotopy-arguments}). Hence a stable solution is asymptotically stable with asymptotic phase (Section~\ref{sec:variational}): nearby solutions converge to a translate of $\tilde{u}$. The same estimate gives a spectral gap (Theorem~\ref{th:spectral-gap}): an explicit lower bound on $|\lambda|$ for every nonzero eigenvalue $\lambda$ of $DF(\tilde{u})$. For a stable solution this gives a lower bound on the rate at which perturbations transverse to the family of translates decay.
The number of negative eigenvalues and the nullspace can in fact already be obtained from the derivative of the augmented system itself, without the new estimate, see Remark~\ref{rem:border-u0}.

The paper is organized as follows. Section~\ref{sec:setting-up-the-system} discusses the variational structure and sets up the augmented system and the contraction argument. It also shows that the unfolding parameters vanish at a fixed point of the contraction map (Section~\ref{sec:beta-zero}). The stability analysis follows in Section~\ref{sec:stability-2D}. Section~\ref{sec:bounds-2D} treats the numerical aspects and derives the bounds needed for the computer-assisted proofs, which are carried out in Section~\ref{sec:results} for four solutions of the SH equation~\eqref{eq:swift-hohenberg-2d}, two stable and two unstable. The computations are done in Julia with IntervalArithmetic.jl~\cite{julia_interval} and RadiiPolynomial.jl~\cite{julia_olivier}, and the code is available at~\cite{julia_blanco_symmetry_stab}. We conclude with an outlook in Section~\ref{sec:conclusion}, and Appendices~\ref{apen:computing-the-bounds} and~\ref{apen:homotopy} contain proofs of technical lemmas.

\section{Setting up the augmented system}
\label{sec:setting-up-the-system}

\subsection{Variational structure}
\label{sec:variational}

The equation $F(u) = 0$ for stationary solutions of~\eqref{eq:swift-hohenberg-2d} is the Euler-Lagrange equation of the energy~\eqref{eq:energy}. We write $(u,v)_2 \bydef \frac{1}{|\om|} \int_{\om} u v^* \, dx$ for the normalized $L^2$ inner product. Since $(I+\Delta)$ is self-adjoint, integration by parts gives
$D\mathcal{E}(u)h = (F(u),h)_2$ for all smooth real-valued periodic functions $u$ and $h$. In other words, $F$ is the $L^2$-gradient of $\mathcal{E}$. Hence the zeros of $F$ are the critical points of $\mathcal{E}$, and $DF(u)$ is the Hessian of $\mathcal{E}$ at $u$, a self-adjoint operator with respect to $(\cdot,\cdot)_2$.

For $s \in \mathbb{R}^2$ let $(\tau_s u)(x) \bydef u(x+s)$ denote translation. The energy is translation invariant, $\mathcal{E}(\tau_s u) = \mathcal{E}(u)$, and differentiating with respect to $s_j$ at $s=0$ gives
\begin{equation}\label{eq:invariance-identity}
  (F(u),\partial_j u)_2 = 0 \qquad \text{for } j=1,2 \quad\text{and all smooth real-valued periodic functions } u .
\end{equation}
We note that~\eqref{eq:invariance-identity} holds for every $u$, not only for solutions.

Translation invariance of $\mathcal{E}$ also explains why solutions are not isolated. Since $\tau_s$ preserves $(\cdot,\cdot)_2$, differentiating $\mathcal{E}(\tau_s u) = \mathcal{E}(u)$ with respect to $u$ gives $F(\tau_s u) = \tau_s F(u)$; this also follows directly from the fact that $F$ does not depend explicitly on $x$. Hence if $F(\tilde{u}) = 0$, then every translate $\tau_s \tilde{u}$ is also a zero. Differentiating $F(\tau_s \tilde{u}) = 0$ with respect to $s_j$ at $s = 0$ gives the following lemma.
\begin{lemma}\label{lem:kernel-of-Df-2d}
Let $\tilde{u}$ be a smooth periodic solution of $F(u) = 0$. Then $\partial_1 \tilde{u}$ and $\partial_2 \tilde{u}$ lie in the nullspace $\mathrm{Null}(DF(\tilde{u}))$ of the Jacobian $DF(\tilde{u})$.
\end{lemma}
In particular, if $\partial_1 \tilde{u}$ and $\partial_2 \tilde{u}$ are linearly independent, then $0$ is an eigenvalue of $DF(\tilde{u})$ of multiplicity at least two. They are linearly dependent if and only if $\alpha_1 \partial_1 \tilde{u} + \alpha_2 \partial_2 \tilde{u} = 0$ for some $\alpha \neq 0$, i.e.~$\tilde{u}$ is constant along the direction $\alpha$ and $\tilde{u}(x) = g(\alpha_2 x_1 - \alpha_1 x_2)$ for a function $g$ of one variable. This is the case only for stripe patterns and constant solutions.

We are interested in the stability of $\tilde{u}$ as a stationary solution of~\eqref{eq:swift-hohenberg-2d}. Since every translate of $\tilde{u}$ is also stationary, stability refers to the family of translates rather than to $\tilde{u}$ alone. If $DF(\tilde{u})$ has no negative eigenvalues and $\mathrm{Null}(DF(\tilde{u})) = \mathrm{span}\{\partial_1 \tilde{u},\partial_2 \tilde{u}\}$, then $\tilde{u}$ is asymptotically stable with asymptotic phase~\cite{henry1981geometric,pruss2009convergence}: solutions of $\partial_t u = -F(u)$ on $\om$ starting close to $\tilde{u}$ converge to $\tau_s \tilde{u}$ for some $s \in \mathbb{R}^2$, which depends on the initial condition. If $DF(\tilde{u})$ has a negative eigenvalue, then $\tilde{u}$ is unstable. The number of negative eigenvalues of $DF(\tilde{u})$, counted with multiplicity, is called the Morse index of $\tilde{u}$.

\subsection{Fourier setting and the augmented system}
\label{sec:fourier-setting}

We write $F$ as
\begin{align}
    F(u) = Lu + G(u), \qquad L \bydef (I + \Delta)^2 + \mu I, \qquad G(u) \bydef -\gamma u^2 + u^3 .\label{def:L-and-G-sh-2d}
\end{align}
We study stationary solutions of~\eqref{eq:swift-hohenberg-2d} with no symmetry enforced, and represent $u$ by its full Fourier series
\[
    u(x) = \sum_{n \in \mathbb{Z}^2} a_n \mathrm{e}^{i \tilde{n} \cdot x}, \qquad \tilde{n} \bydef \left(\frac{n_1\pi}{d_1},\frac{n_2\pi}{d_2}\right).
\]
We identify $u$ with its sequence of Fourier coefficients $(a_n)_{n \in \mathbb{Z}^2}$. We work with complex sequences; real-valued functions correspond to conjugate symmetric sequences, i.e.~$a_{-n} = a_n^*$, where $*$ denotes complex conjugation. In these coordinates $L$ is diagonal,
\[
    (Lu)_n = L_{nn} a_n, \qquad L_{nn} \bydef (1 - |\tilde{n}|^2)^2 + \mu.
\]
The nonlinearity $G$ is given by discrete convolutions. Under this identification the inner product of Section~\ref{sec:variational} becomes the inner product of $\ell^2(\mathbb{Z}^2)$:
\[
    (u,v)_2 = \sum_{n \in \mathbb{Z}^2} a_n b_n^* \qquad \text{for } u = (a_n)_{n \in \mathbb{Z}^2}, ~ v = (b_n)_{n \in \mathbb{Z}^2} .
\]
We write $\|u\|_2 \bydef (u,u)_2^{1/2} = \bigl(\sum_{n \in \mathbb{Z}^2} |a_n|^2\bigr)^{1/2}$ and $\|u\|_1 \bydef \sum_{n \in \mathbb{Z}^2} |a_n|$. We have $\|u\|_2 \leq \|u\|_1$.

Since solutions are not isolated, we augment $F$ with two unfolding parameters $\beta_1$ and $\beta_2$ and two phase conditions. Let $u_0$ be a real-valued \emph{approximate} solution, i.e.~$F(u_0) \approx 0$. For $w \bydef (\beta_1,\beta_2,u)$ we consider the augmented system
\begin{align}
    \overline{F}(w) \bydef \begin{bmatrix}
        (u,\partial_1 u_0)_2 \\
        (u,\partial_2 u_0)_2 \\
        F(u) + \beta_1 \partial_1 u_0 + \beta_2 \partial_2 u_0
    \end{bmatrix}.\label{def:F_bar-2d}
\end{align}

On $\mathbb{C}^2 \times \ell^2(\mathbb{Z}^2)$ we use the inner product 
\begin{align}
\label{e:inner2}
((\beta,u),(\gamma,v)) \bydef \beta_1\gamma_1^* + \beta_2\gamma_2^* + (u,v)_2.
\end{align}
The augmented system inherits the variational structure of Section~\ref{sec:variational}. With
\[
  \overline{\mathcal{E}}(\beta_1,\beta_2,u) \bydef \mathcal{E}(u) + \beta_1 (u,\partial_1 u_0)_2 + \beta_2 (u,\partial_2 u_0)_2
\]
we have $\overline{F} = \nabla \overline{\mathcal{E}}$ for real $\beta$ and real-valued $u$. The unfolding parameters are thus Lagrange multipliers for the phase conditions $(u,\partial_1 u_0)_2 = (u,\partial_2 u_0)_2 = 0$, which fix the position of $u$ along the family of translates. In particular
\[
  D\overline{F}(w) = \begin{bmatrix}
        0 & 0 & (\partial_1 u_0)^* \\
        0 & 0 & (\partial_2 u_0)^* \\
        \partial_1 u_0 & \partial_2 u_0 & DF(u)
    \end{bmatrix}
\]
is self-adjoint for real-valued $u$, where $(\partial_j u_0)^*$ denotes the functional $v \mapsto (v,\partial_j u_0)_2$. The point $w_0 = (0,0,u_0)$ satisfies the phase conditions exactly, since $(u_0,\partial_j u_0)_2 = \frac{1}{2|\om|}\int_{\om} \partial_j (u_0^2) \, dx = 0$.

For $\nu = (\nu_1,\nu_2) > 1$, where here and below inequalities between pairs are meant componentwise, let $\ell^1_\nu$ be the space of sequences $u = (a_n)_{n \in \mathbb{Z}^2}$ with
\[
    \|u\|_{\ell^1_\nu} \bydef \sum_{n \in \mathbb{Z}^2} |a_n| \nu_1^{|n_1|}\nu_2^{|n_2|} < \infty .
\]
For the augmented system we use the two components
\[
    X_1 \bydef \mathbb{C}^2 \quad \text{with } \|\beta\|_{X_1} \bydef |\beta_1| + |\beta_2|, \qquad X_2 \bydef \ell^1_\nu \quad \text{with } \|\cdot\|_{X_2} \bydef \|\cdot\|_{\ell^1_\nu}.
\]
Since $L$ is of fourth order, the range of $\overline{F}$ is less regular. We let $Y_1 \bydef \mathbb{C}^2$ and
\[
    Y_2 \bydef \Bigl\{v = (v_{n})_{n \in \mathbb{Z}^2}, \|v\|_{Y_2} \bydef \sum_{n \in \mathbb{Z}^2} |v_n| \frac{\nu_1^{|n_1|}\nu_2^{|n_2|}}{(1+|\tilde{n}|^2)^2} < \infty\Bigr\}.
\]
Since $\ell^1_\nu$ is a Banach algebra under discrete convolution, $\overline{F}$ is a smooth map from $X_1 \times X_2$ to $Y_1 \times Y_2$.
We denote by $P_1$ and $P_2$ the projections of $X_1 \times X_2$ onto $X_1$ and $X_2$, and also those of $Y_1 \times Y_2$ onto $Y_1$ and $Y_2$. Bounded linear operators $\Gamma$ between such product spaces are written as $2 \times 2$ block matrices, where the $(i,j)$ block $P_i \Gamma P_j$ maps the $j$-th component to the $i$-th component; for instance, for $\Gamma \in \mathcal{B}(Y_1 \times Y_2, X_1 \times X_2)$ we have $P_i \Gamma P_j \in \mathcal{B}(Y_j,X_i)$.

\subsection{The radii polynomial theorem}
\label{sec:radii-polynomial}
Let $u_0$ be the real-valued approximate solution of Section~\ref{sec:fourier-setting}. In practice it is a trigonometric polynomial computed numerically, see Section~\ref{sec:numerical-aspects}. We expect the unfolding parameters to vanish at a solution (this is shown in Section~\ref{sec:beta-zero}), and $u_0$ satisfies the phase conditions exactly, so $w_0 \bydef (0,0,u_0)$ is an approximate zero of $\overline{F}$. 

For the approximate inverse of $D\overline{F}(w_0)$ we split the Fourier modes at a truncation level $N = (N_1,N_2)$. Let $I^N \bydef \{n \in \mathbb{Z}^2 : |n_1| \leq N_1, |n_2| \leq N_2\}$, and let $\pi^{\le N}$ and $\pi^{>N}$ be the projections that keep the Fourier coefficients with $n \in I^N$ and $n \notin I^N$, respectively. 
Sums and multiples of truncation levels, such as $2N$ or $N + N_0$, are taken componentwise.
We denote by $N_0$ the truncation level of $u_0$, i.e.~$u_0$ has nonzero Fourier coefficients only for $n \in I^{N_0}$. We choose $N$ so large that $L_{nn} > 0$ for all $n \notin I^N$ (for $\mu>0$ this is trivial). The approximate inverse $\overline{A} : Y_1 \times Y_2 \to X_1 \times X_2$ is chosen to be block diagonal,
\[
    \overline{A} = \overline{A}^N \oplus \pi^{>N} L^{-1} .
\]
Here $\overline{A}^N$ is a self-adjoint matrix acting on $\beta$ and on the Fourier coefficients with $n \in I^N$, which approximates the inverse of the corresponding finite block of $D\overline{F}(w_0)$. On the remaining Fourier coefficients $\overline{A}$ acts as the diagonal operator $\pi^{>N} L^{-1}$, which has positive entries $1/L_{nn}$. In particular $\overline{A}$ is self-adjoint with respect to the inner product of Section~\ref{sec:fourier-setting}. Further details are given in Section~\ref{sec:numerical-aspects}. We introduce the fixed point map
\begin{align}
    T(w) \bydef w - \overline{A}\,\overline{F}(w).
\end{align}
For $r \bydef (r_1,r_2) > 0$ we define the box
\begin{align}
    B_r(w_0) \bydef \{w \in X_1 \times X_2 : \|P_1 (w - w_0)\|_{X_1} \leq r_1, ~ \|P_2 (w - w_0)\|_{X_2} \leq r_2\}.
\end{align}
We wish to prove that $T$ maps $B_r(w_0)$ to itself and that $T$ is a contraction. For $\chi = (\chi_1,\chi_2) > 0$ we use the weighted norm
\[
    \|h\|_\chi \bydef \chi_1 \|P_1 h\|_{X_1} + \chi_2 \|P_2 h\|_{X_2}, \qquad h \in X_1 \times X_2 ,
\]
on $X_1 \times X_2$. The weights $\chi$ are chosen such that $T$ is a contraction in this norm, see~\eqref{condition-contraction} in the following theorem, and a norm of the same form, with other weights $\xi$, is used in Section~\ref{sec:stability-2D}. The theorem is a version of the radii polynomial approach, of which many variations have been used, see e.g.~\cite{day_lessard_mischaikow_2007,vandenberg_introduction_2018,JB_symmetries_2}. We include a proof for completeness.
\begin{theorem}\label{th:radii-polynomial}
Let $\overline{A} : Y_1 \times Y_2 \to X_1 \times X_2$ be the approximate inverse above and let $r^* = (r_1^*,r_2^*) > 0$. Let $Y^{[i]}, Z_j^{[i]}, W_{jk}^{[i]} \geq 0$ for $i,j,k = 1,2$ be such that
\begin{equation}\label{eq:definition-Y0-Z1-Z2}
\begin{aligned}
    \|P_i (T(w_0) - w_0)\|_{X_i} &\leq Y^{[i]}, \\
    \|P_i DT(w_0) P_j\|_{\mathcal{B}(X_j,X_i)} &\leq Z_j^{[i]}, \\
    \|P_i\bigl(DT(w) - DT(w_0)\bigr)P_j\|_{\mathcal{B}(X_j,X_i)} &\leq W_{j1}^{[i]} \hat{r}_1 + W_{j2}^{[i]} \hat{r}_2 \qquad \text{for all } w \in B_{\hat{r}}(w_0) \text{ and all } 0 < \hat{r} \leq r^*.
\end{aligned}
\end{equation}
Let $r = (r_1,r_2)$ with $0 < r \leq r^*$. If
\begin{equation}\label{condition-radii-polynomial}
    Y^{[i]} + Z_1^{[i]} r_1 + Z_2^{[i]} r_2 + \frac{1}{2} \sum_{j,k=1}^2 W_{jk}^{[i]} r_j r_k < r_i \qquad \text{for } i = 1,2,
\end{equation}
and there exists $\chi = (\chi_1,\chi_2) > 0$ such that
\begin{equation}\label{condition-contraction}
    \chi_1 H_j^{[1]}(r) + \chi_2 H_j^{[2]}(r) < \chi_j \qquad \text{for } j = 1,2,
\end{equation}
where
\begin{equation}\label{eq:definition-U}
    H_j^{[i]}(r) \bydef Z_j^{[i]} + W_{j1}^{[i]} r_1 + W_{j2}^{[i]} r_2 , \qquad i,j = 1,2,
\end{equation}
then there exists a unique $\tilde{w} \in B_r(w_0)$ such that $\overline{F}(\tilde{w})=0$.
\end{theorem}
\begin{proof}
Let $w \in B_r(w_0)$, let $h \bydef w - w_0$ and denote $\hat{h}_k \bydef \|P_k h\|_{X_k} \leq r_k$. Since
\[
    T(w) - w_0 = T(w_0) - w_0 + DT(w_0)h + \int_0^1 \bigl(DT(w_0 + th) - DT(w_0)\bigr) h \, dt
\]
and $w_0 + th \in B_{tr}(w_0)$ for $t \in (0,1]$, we obtain
\[
  \|P_i(T(w) - w_0)\|_{X_i} \leq Y^{[i]} + \sum_{j=1}^2 Z_j^{[i]} \hat{h}_j + \sum_{j,k=1}^2  W_{jk}^{[i]} \int_0^1 t\, r_k \hat{h}_j  \, dt  \leq Y^{[i]} + \sum_{j=1}^2 Z_j^{[i]} r_j + \frac{1}{2}\sum_{j,k=1}^2 W_{jk}^{[i]} r_j r_k < r_i .
\]
Hence $T$ maps $B_r(w_0)$ into itself. 

For $w \in B_r(w_0)$ we have $\|P_i DT(w) P_j\|_{\mathcal{B}(X_j,X_i)} \leq H_j^{[i]}(r)$, so for any $h \in X_1 \times X_2$
\[
  \|DT(w)h\|_\chi \leq \sum_{i=1}^2 \chi_i \sum_{j=1}^2 H_j^{[i]}(r) \|P_j h\|_{X_j} = \sum_{j=1}^2 \bigl(\chi_1 H_j^{[1]}(r) + \chi_2 H_j^{[2]}(r)\bigr) \|P_j h\|_{X_j} \leq \theta \|h\|_\chi ,
\]
where $\theta \bydef \max_{j} (\chi_1 H_j^{[1]}(r) + \chi_2 H_j^{[2]}(r)) / \chi_j < 1$. Since $B_r(w_0)$ is convex and closed, $T$ is a contraction on $B_r(w_0)$ with respect to $\|\cdot\|_\chi$, and it has a unique fixed point $\tilde{w} \in B_r(w_0)$. The bound $\|I - \overline{A}D\overline{F}(w_0)\|_\chi \leq \theta < 1$ and the block structure of $\overline{A}$ imply that $\overline{A}$ is injective, so the zeros of $\overline{F}$ in $B_r(w_0)$ are exactly the fixed points of $T$.
\end{proof}

Since $\overline{F}$ is affine in $\beta$, the derivative $D\overline{F}(w)$ depends on $w$ only through $u$, and its first block column does not depend on $w$ at all. Hence $W_{jk}^{[i]} = 0$ unless $j = k = 2$, and only $W_{22}^{[1]}$ and $W_{22}^{[2]}$ need to be computed (Lemma~\ref{lem:Z2}).

\begin{remark}
Since $DT(\tilde{w}) = I - \overline{A}D\overline{F}(\tilde{w})$ and $\tilde{w} \in B_r(w_0)$, the proof of Theorem~\ref{th:radii-polynomial} also gives
\begin{equation}\label{second-assumption-in-homotopy-lemma}
    \|I - \overline{A}D\overline{F}(\tilde{w})\|_\chi < 1 .
\end{equation}
This bound is used in Remark~\ref{rem:border-u0}.
\end{remark}

\begin{remark}\label{rem:contraction-condition}
Let $r$ be fixed such that~\eqref{condition-radii-polynomial} holds, and let $\mathcal{H} = \mathcal{H}(r)$ be the $2 \times 2$ matrix with entries $\mathcal{H}_{ij} \bydef H_j^{[i]}(r)$. Then condition~\eqref{condition-contraction} states that $\chi^T \mathcal{H} < \chi^T$ componentwise. Often the condition $\mathcal{H}\eta < \eta$ for some $\eta > 0$, i.e.,
\[
    H_1^{[i]}(r) \eta_1 + H_2^{[i]}(r) \eta_2 < \eta_i \qquad \text{for } i = 1,2,
\]
is used instead, see e.g.~\cite{JB_symmetries_2,van2021spontaneous}; it corresponds to the weighted maximum norm $\max_i \|P_i h\|_{X_i}/\eta_i$. We use~\eqref{condition-contraction}, which corresponds to the weighted $\ell^1$ norm $\|\cdot\|_\chi$, because we need a norm of this form later, in Lemma~\ref{lem:schur} and Theorem~\ref{th:spectral-gap}. 

Since $\mathcal{H}$ has nonnegative entries, both conditions are equivalent to $\rho(\mathcal{H}) < 1$, where the spectral radius $\rho(\mathcal{H})$ is the largest absolute value of the eigenvalues of $\mathcal{H}$; they are instances of the Collatz-Wielandt formula for the spectral radius of a nonnegative matrix, applied to $\mathcal{H}^T$ and $\mathcal{H}$, respectively, see e.g.~\cite[Chapter~8]{horn_johnson_2013}. If $\rho(\mathcal{H}) < 1$, then $\chi \bydef (I - \mathcal{H}^T)^{-1}(1,1)^T$ has positive entries and satisfies $\chi^T \mathcal{H} = \chi^T - (1,1)$; in practice we compute this $\chi$ in floating point arithmetic and verify~\eqref{condition-contraction} with interval arithmetic.\end{remark}

Theorem~\ref{th:radii-polynomial} is formulated for complex sequences. For $u = (a_n)_{n \in \mathbb{Z}^2}$ let $Su \bydef (a_{-n}^*)_{n \in \mathbb{Z}^2}$ and $S(\beta,u) \bydef (\beta^*,Su)$. The conjugate symmetric sequences, which correspond to real-valued functions, are then the fixed points of $S$. Since $L_{nn}$ and the weights are real and even in $n$, $S$ preserves the norms and $F(Su) = SF(u)$. If $u_0$ is conjugate symmetric, then also $\overline{F}(Sw) = S\overline{F}(w)$.
\begin{lemma}\label{lem:real-solution}
Let the assumptions of Theorem~\ref{th:radii-polynomial} hold and let $u_0$ be conjugate symmetric. Then the zero $\tilde{w} = (\tilde{\beta},\tilde{u})$ of $\overline{F}$ in $B_r(w_0)$ satisfies $\tilde{\beta} \in \mathbb{R}^2$, and $\tilde{u}$ is conjugate symmetric, i.e.~real-valued.
\end{lemma}
\begin{proof}
Since $S$ maps $B_r(w_0)$ to itself, $S\tilde{w}$ is also a zero of $\overline{F}$ in $B_r(w_0)$, so $S\tilde{w} = \tilde{w}$ by uniqueness.
\end{proof}
If the computed $u_0$ is conjugate symmetric only up to a small defect, we apply the theorem to $\frac{1}{2}(u_0 + Su_0)$.

\subsection{Vanishing of the unfolding parameters}
\label{sec:beta-zero}

Theorem~\ref{th:radii-polynomial} provides a zero $\tilde{w} = (\tilde{\beta}_1,\tilde{\beta}_2,\tilde{u})$ of the augmented system~\eqref{def:F_bar-2d}. It remains to show that $\tilde{\beta} = 0$, so that $F(\tilde{u}) = 0$. We prove that this follows from the translation invariance, provided the matrix
\begin{equation}\label{eq:def-C}
  M \bydef \begin{bmatrix}
    (\partial_1 u_0, \partial_1 \tilde{u})_2 & (\partial_2 u_0, \partial_1 \tilde{u})_2 \\
    (\partial_1 u_0, \partial_2 \tilde{u})_2 & (\partial_2 u_0, \partial_2 \tilde{u})_2
  \end{bmatrix}
\end{equation}
is invertible.
\begin{theorem}\label{beta=0-2d}
Let $u_0$ be conjugate symmetric and let $\tilde{w} = (\tilde{\beta}_1,\tilde{\beta}_2,\tilde{u}) \in B_r(w_0)$ be the zero of $\overline{F}$ provided by Theorem~\ref{th:radii-polynomial}. If $M$ is invertible, then $\tilde{\beta} = 0$ and $\tilde{u}$ is a real-valued zero of $F$.
\end{theorem}
\begin{proof}
By Lemma~\ref{lem:real-solution}, $\tilde{u}$ is real-valued; since $\nu > 1$, $\tilde{u}$ is also smooth. The last component of $\overline{F}(\tilde{w}) = 0$ reads $F(\tilde{u}) + \tilde{\beta}_1 \partial_1 u_0 + \tilde{\beta}_2 \partial_2 u_0 = 0$. Taking the inner product with $\partial_i \tilde{u}$ and using 
\eqref{eq:invariance-identity}, we obtain
\[
  0 = (F(\tilde{u}),\partial_i \tilde{u})_2 + \sum_{j=1}^2 \tilde{\beta}_j (\partial_j u_0, \partial_i \tilde{u})_2 = (M\tilde{\beta})_i , \qquad i=1,2.
\]
Since $M$ is invertible by assumption, $\tilde{\beta} = 0$, so the last component of $\overline{F}(\tilde{w}) = 0$ reduces to $F(\tilde{u}) = 0$.
\end{proof}

The invertibility of $M$ also shows that $\tilde{u}$ is locally unique up to translations. For any $u$ close to $\tilde{u}$ consider the map $s \mapsto \bigl((\tau_s u,\partial_1 u_0)_2, (\tau_s u,\partial_2 u_0)_2\bigr)$, which is smooth, since $(\tau_s u,\partial_j u_0)_2 = (u,\tau_{-s}\partial_j u_0)_2$ and $u_0$ is smooth. At $s=0$ and $u=\tilde{u}$ it vanishes and its derivative is $M^T$, so by the implicit function theorem there is a small $s$ such that $\tau_s u$ satisfies the phase conditions. If moreover $F(u) = 0$, then also $F(\tau_s u) = 0$, so $(0,0,\tau_s u)$ is a zero of $\overline{F}$ in $B_r(w_0)$. By uniqueness, $\tau_s u = \tilde{u}$, that is, $u = \tau_{-s}\tilde{u}$.

The matrix $M$ involves the unknown solution $\tilde{u}$. We write $\tilde{u} = u_0 + v$ with $v \bydef P_2(\tilde{w} - w_0)$, so that $\|v\|_{\ell^1_\nu} \leq r_2$ since $\tilde{w} \in B_r(w_0)$. Let $M = M_0 + M_\varepsilon$ with
\[
  (M_0)_{ij} \bydef (\partial_j u_0, \partial_i u_0)_2 , \qquad (M_\varepsilon)_{ij} \bydef (\partial_j u_0, \partial_i v)_2 .
\]
The matrix $M_0$ involves only the approximate solution $u_0$, and not $\tilde{u}$, and it is positive definite if $\partial_1 u_0$ and $\partial_2 u_0$ are linearly independent, since $\alpha^T M_0 \alpha = \|\alpha_1 \partial_1 u_0 + \alpha_2 \partial_2 u_0\|_2^2$. By Section~\ref{sec:variational} this is the case unless $u_0$ is a stripe pattern or constant. The perturbation $M_\varepsilon$ is controlled by the following estimate.
\begin{lemma}\label{lem:derivative-estimate}
For $v \in \ell^1_\nu$ and $i=1,2$, we have 
$\|\partial_i v\|_1 \leq \kappa \|v\|_{\ell^1_\nu}$, where  
\begin{align*}
  \kappa \bydef \max_{k=1,2} \frac{\pi}{\mathrm{e}\, d_k \ln \nu_k} .
\end{align*}
\end{lemma}
\begin{proof}
We have $\|\partial_i v\|_1 = \sum_{n} |\tilde{n}_i| |v_n| \leq \sup_{n} |\tilde{n}_i| \nu_1^{-|n_1|}\nu_2^{-|n_2|} \, \|v\|_{\ell^1_\nu}$ and
\[
  |\tilde{n}_i| \nu_1^{-|n_1|}\nu_2^{-|n_2|} \leq \frac{\pi}{d_i} |n_i| \nu_i^{-|n_i|} \leq \frac{\pi}{d_i} \max_{t \geq 0} t \nu_i^{-t} = \frac{\pi}{\mathrm{e}\, d_i \ln \nu_i} \leq \kappa . \qedhere
\]
\end{proof}
By using $\|\cdot\|_2 \leq \|\cdot\|_1$, we see that $|(M_\varepsilon)_{ij}| \leq \|\partial_j u_0\|_2 \|\partial_i v\|_1 \leq \kappa \|\partial_j u_0\|_2 \, r_2$. Hence $M$ lies in the interval matrix
\begin{equation}\label{eq:C-enclosure}
  \mathbf{M} \bydef M_0 + \kappa r_2 \begin{bmatrix}
     \|\partial_1 u_0\|_2 & \|\partial_2 u_0\|_2 \\
     \|\partial_1 u_0\|_2 & \|\partial_2 u_0\|_2
    \end{bmatrix} [-1,1],
\end{equation}
where the interval $[-1,1]$ multiplies each entry. Since $u_0$ is a trigonometric polynomial, the entries of $M_0$ and the norms $\|\partial_j u_0\|_2$ are finite sums, which we enclose with interval arithmetic. To apply Theorem~\ref{beta=0-2d} we verify with interval arithmetic that $0 \notin \det \mathbf{M}$. Since $M_0$ is positive definite and $\kappa r_2$ is small in practice, this is a mild condition. Alternatively, to limit overestimation in interval arithmetic, one may compute $R \approx M_0^{-1}$ in floating point and verify with interval arithmetic that $\|I - R\,\mathbf{M}\| < 1$ for some induced matrix norm, implying that $M$ is invertible.

One could also integrate by parts, $(M_\varepsilon)_{ij} = -(\partial_i\partial_j u_0, v)_2$, and bound $M_\varepsilon$ by a finite computation without $\kappa$, since $u_0$ has finitely many nonzero Fourier coefficients. We do not use this, since $\kappa$ is needed anyway for $Z_{Q,2}^{[i]}$ in Section~\ref{sec:Zq}.

\section{Stability}
\label{sec:stability-2D}

Let $\tilde{w} = (0,0,\tilde{u})$ be the zero of $\overline{F}$ obtained from Theorems~\ref{th:radii-polynomial} and~\ref{beta=0-2d}. Stability of $\tilde{u}$ is decided by its Morse index $\mathcal{N}(DF(\tilde{u}))$, the number of negative eigenvalues of $DF(\tilde{u})$, and by the nullspace of $DF(\tilde{u})$.

The main result of this section is Corollary~\ref{lem:homotopy-arguments}. It reduces the number of unstable directions to the number of negative eigenvalues of a finite matrix, which is determined rigorously in Section~\ref{sec:numerical-aspects}. Along the way we show that the only neutral directions of $DF(\tilde{u})$ are the translations (Lemma~\ref{lem:eigenvalues-of-Q-2d}). We proceed in the following steps.
\begin{enumerate}
\item Since $DF(\tilde{u})$ is singular, its negative eigenvalues cannot be counted with an approximate inverse directly. We therefore extend $DF(\tilde{u})$ by two rows and two columns, as in the augmented system, but with the exact (rather than approximate) translation directions in the border. If the translation directions are linearly independent and span the nullspace of $DF(\tilde{u})$, this extended operator $Q$ has the same eigenvalues as $DF(\tilde{u})$, except that the double eigenvalue zero is replaced by two positive and two negative eigenvalues (Lemma~\ref{lem:bordered}).
\item The approximate inverse from the existence proof is also an approximate inverse of the extended operator $Q$, provided one additional inequality holds that the computer can check (Theorem~\ref{th:radii-polynomial-q}). Together with the invertibility of $M$ and the fact that the relevant eigenvectors decay rapidly (Lemma~\ref{lem:regularity}), this verifies these two conditions (Lemma~\ref{lem:eigenvalues-of-Q-2d}).
\item The extended operator $Q$ has as many negative eigenvalues as the approximate inverse (Proposition~\ref{prop:homotopy}). Since the tail of the approximate inverse is positive, these are the negative eigenvalues of a finite matrix, and together with step~1 this gives the count for $DF(\tilde{u})$ (Corollary~\ref{lem:homotopy-arguments}). Moreover, the approximate inverse is close to the inverse of $Q$, so its norm gives a lower bound on the absolute values of the nonzero eigenvalues of $DF(\tilde{u})$ (Theorem~\ref{th:spectral-gap}).
\item The count and the statement on the neutral directions can also be obtained directly from the derivative of the augmented system, without the additional inequality of step~2 (Remark~\ref{rem:border-u0}).
\end{enumerate}

Stability is a property of $DF(\tilde{u})$ as the Hessian of $\mathcal{E}$ with respect to the $L^2$ inner product, see Section~\ref{sec:variational}, and the tools we need (self-adjointness, orthogonal decompositions, and the characterization of $\mathcal{N}$ by quadratic forms) all rely on this inner product. In this section we therefore work in $\ell^2 = \ell^2(\mathbb{Z}^2)$ with the inner product $(\cdot,\cdot)_2$, and in $\mathbb{C}^2 \times \ell^2$ with the inner product~\eqref{e:inner2}. From now on we write $\|\cdot\|_2$ also for the norm of this inner product, and $\|\cdot\|_{\mathcal{B}(\ell^2)}$ also for the corresponding operator norm on $\mathbb{C}^2 \times \ell^2$. The linear operators whose eigenvalues we study below are self-adjoint on these spaces. We work over $\mathbb{C}$, as in the computations. Since $\tilde{u}$ is real-valued, $DF(\tilde{u})$ commutes with the conjugation $S$ of Section~\ref{sec:radii-polynomial}, whose fixed points are the conjugate symmetric sequences, which correspond to real-valued functions. Hence every complex eigenspace of $DF(\tilde{u})$ is invariant under $S$ and spanned by conjugate symmetric eigenvectors, so the eigenvalues and their multiplicities are the same over $\mathbb{C}$ as for real-valued functions. Since $L_{nn} \to \infty$ and $DF(\tilde{u}) - L$ is bounded, $DF(\tilde{u})$ has compact resolvent and finitely many negative eigenvalues.

Every translate of $\tilde{u}$ is again a stationary solution, hence the translation directions $\partial_1\tilde{u}$ and $\partial_2\tilde{u}$ lie in the nullspace of $DF(\tilde{u})$, see Lemma~\ref{lem:kernel-of-Df-2d}. As already mentioned, $DF(\tilde{u})$ is not invertible, and its negative eigenvalues cannot be counted directly with an approximate inverse as in~\cite{JB_symmetries_2}. Instead, for $w = (\beta,u)$ we consider the extended, bordered operator
\begin{equation}\label{def:Q-2d}
    Q(w) \bydef \begin{bmatrix}
        0 & 0 & (\partial_1 u)^* \\
        0 & 0 & (\partial_2 u)^* \\
        \partial_1 u & \partial_2 u & DF(u)
    \end{bmatrix},
\end{equation}
which differs from $D\overline{F}(w)$ only in the border: $\partial_j u$ replaces $\partial_j u_0$. At $w = \tilde{w}$ the border consists of the null vectors $\partial_1 \tilde{u}$, $\partial_2 \tilde{u}$ of $DF(\tilde{u})$, and this makes the spectrum of $Q(\tilde{w})$ explicit in terms of that of $DF(\tilde{u})$.

Let $K \bydef \mathrm{span}\{\partial_1\tilde{u},\partial_2\tilde{u}\}$ be the space of translation directions; we will show that $K$ is the whole nullspace of $DF(\tilde{u})$. The four eigenvalues of $Q(\tilde{w})$ that replace the eigenvalue zero on $K$ are determined by the matrix $\Phi$ with entries $\Phi_{ij} \bydef (\partial_j\tilde{u},\partial_i\tilde{u})_2$. It is the analogue, for the exact border, of the matrix $M$ of Section~\ref{sec:beta-zero}. If $\partial_1\tilde{u}$ and $\partial_2\tilde{u}$ are linearly independent, then $\Phi$ is positive definite. We denote its eigenvalues by $\lambda_1^2,\lambda_2^2 > 0$.
\begin{lemma}\label{lem:bordered}
Let $\partial_1\tilde{u}$ and $\partial_2\tilde{u}$ be linearly independent. Then the eigenvalues of $Q(\tilde{w})$ on $\mathbb{C}^2 \times \ell^2$ are $\pm\lambda_1$, $\pm\lambda_2$ and the eigenvalues of $DF(\tilde{u})$ restricted to $K^\perp$. The latter include all nonzero eigenvalues of $DF(\tilde{u})$. If $Q(\tilde{w})$ is injective on $\mathbb{C}^2 \times \ell^2$, then $\mathrm{Null}(DF(\tilde{u})) = K$ and $\mathcal{N}(DF(\tilde{u})) = \mathcal{N}(Q(\tilde{w})) - 2$.
\end{lemma}
\begin{proof}
Recall that $DF(\tilde{u})\partial_j\tilde{u} = 0$. The idea is that the exact border decouples: $\mathbb{C}^2 \times \ell^2$ is the orthogonal sum of the four-dimensional space $\mathbb{C}^2 \times K$, which consists of the two unfolding directions and the two translation directions, and $\{0\} \times K^\perp$. We show that $Q(\tilde{w})$ maps each of these spaces into itself; its spectrum is then the union of the spectra on the two parts.

On $\{0\} \times K^\perp$: for $y \in K^\perp$ the first two components of $Q(\tilde{w})(0,y)$, namely $(y,\partial_1\tilde{u})_2$ and $(y,\partial_2\tilde{u})_2$, vanish, so $Q(\tilde{w})(0,y) = (0,DF(\tilde{u})y)$, and $DF(\tilde{u})y \in K^\perp$ since $(DF(\tilde{u})y,\partial_j\tilde{u})_2 = (y,DF(\tilde{u})\partial_j\tilde{u})_2 = 0$. Hence $Q(\tilde{w})$ acts on $\{0\} \times K^\perp$ as $DF(\tilde{u})$ acts on $K^\perp$.

On $\mathbb{C}^2 \times K$: here $Q(\tilde{w})$ acts as 
\[
  Q(\tilde{w})(\beta, \alpha_1\partial_1\tilde{u} + \alpha_2\partial_2\tilde{u}) = (\Phi\alpha, \beta_1\partial_1\tilde{u} + \beta_2\partial_2\tilde{u}) \qquad \text{for all } \alpha,\beta \in \mathbb{C}^2 ,
\]
since $DF(\tilde{u})\partial_j\tilde{u} = 0$ and $(\alpha_1\partial_1\tilde{u} + \alpha_2\partial_2\tilde{u},\partial_i\tilde{u})_2 = (\Phi\alpha)_i$. This means that in these coordinates the action of $Q(\tilde{w})$ on $\mathbb{C}^2 \times K$ is described by the $4 \times 4$ matrix with blocks $\Phi$ and $I$ off the diagonal. If $\Phi v = \lambda^2 v$ with $\lambda > 0$, then $(\pm\lambda v, v)$ is an eigenvector with eigenvalue $\pm\lambda$, so the eigenvalues on this part are $\pm\lambda_1,\pm\lambda_2$.

Finally, $DF(\tilde{u})$ is the orthogonal sum of $0$ on $K$ and its restriction to $K^\perp$. This gives the remaining statements, since exactly two of $\pm\lambda_1,\pm\lambda_2$ are negative.
\end{proof}

The operator $Q(\tilde{w})$ is not known explicitly, while the computer-assisted proof controls $\overline{A}$ and $D\overline{F}$. The following extension of Theorem~\ref{th:radii-polynomial} relates $\overline{A}$ to $Q(\tilde{w})$, in the weighted norm $\|\cdot\|_\xi$, defined as $\|\cdot\|_\chi$ in Section~\ref{sec:radii-polynomial} with weights $\xi$ in place of $\chi$.
\begin{theorem}\label{th:radii-polynomial-q}
Let $\overline{A}$, the bounds~\eqref{eq:definition-Y0-Z1-Z2} and $r$ be as in Theorem~\ref{th:radii-polynomial}, where $r$ satisfies~\eqref{condition-radii-polynomial}. Let $Z_{Q,j}^{[i]}(r) \geq 0$ for $i,j = 1,2$ be such that
\begin{equation}\label{eq:definition-Zq}
    \|P_i\overline{A}(Q(w) - D\overline{F}(w))P_j\|_{\mathcal{B}(X_j,X_i)} \leq Z_{Q,j}^{[i]}(r) \qquad \text{for all } w \in B_r(w_0).
\end{equation}
If there exists $\xi = (\xi_1,\xi_2) > 0$ such that
\begin{equation}\label{condition-radii-polynomial-q}
    \xi_1 \bigl(H_j^{[1]}(r) + Z_{Q,j}^{[1]}(r)\bigr) + \xi_2 \bigl(H_j^{[2]}(r) + Z_{Q,j}^{[2]}(r)\bigr) < \xi_j \qquad \text{for } j = 1,2,
\end{equation}
with $H_j^{[i]}(r)$ as in~\eqref{eq:definition-U}, then $\overline{F}$ has a unique zero $\tilde{w} \in B_r(w_0)$, and
\begin{equation}
    \|I - \overline{A}Q(\tilde{w})\|_\xi \leq \delta \bydef \max_{j} \frac{\xi_1 \bigl(H_j^{[1]}(r) + Z_{Q,j}^{[1]}(r)\bigr) + \xi_2 \bigl(H_j^{[2]}(r) + Z_{Q,j}^{[2]}(r)\bigr)}{\xi_j} < 1. \label{first-assumption-in-homotopy-lemma}
\end{equation}
\end{theorem}
\begin{proof}
Since $Z_{Q,j}^{[i]}(r) \geq 0$, condition~\eqref{condition-radii-polynomial-q} implies~\eqref{condition-contraction} with $\chi = \xi$, so Theorem~\ref{th:radii-polynomial} applies. Since $I - \overline{A}D\overline{F}(\tilde{w}) = DT(\tilde{w})$ and $\tilde{w} \in B_r(w_0)$, the blocks $P_i DT(\tilde{w}) P_j$ are bounded by $H_j^{[i]}(r)$, as shown in the proof of Theorem~\ref{th:radii-polynomial}. By~\eqref{eq:definition-Zq}, the blocks of $\overline{A}(Q(\tilde{w}) - D\overline{F}(\tilde{w}))$ are bounded by $Z_{Q,j}^{[i]}(r)$, so the blocks of $I - \overline{A}Q(\tilde{w}) = DT(\tilde{w}) - \overline{A}(Q(\tilde{w}) - D\overline{F}(\tilde{w}))$ are bounded by $H_j^{[i]}(r) + Z_{Q,j}^{[i]}(r)$, and the computation in the proof of Theorem~\ref{th:radii-polynomial} gives~\eqref{first-assumption-in-homotopy-lemma}.
\end{proof}

The weight vector $\xi$ in Theorem~\ref{th:radii-polynomial-q} is in general not the vector $\chi$ of Remark~\ref{rem:contraction-condition}: it is computed in the same way, with $H_j^{[i]}(r) + Z_{Q,j}^{[i]}(r)$ in place of $H_j^{[i]}(r)$.

The estimate~\eqref{first-assumption-in-homotopy-lemma} says that $\overline{A}$ is an approximate inverse not only of $D\overline{F}(w_0)$, for which it was constructed, but also of $Q(\tilde{w})$ at the solution itself. It gives the injectivity of $Q(\tilde{w})$ on $X_1 \times X_2$ required in Lemma~\ref{lem:eigenvalues-of-Q-2d} below, it allows us to count the negative eigenvalues of $Q(\tilde{w})$ through those of $\overline{A}$ (Proposition~\ref{prop:homotopy}), and the number $\delta$ enters the spectral gap (Theorem~\ref{th:spectral-gap}). The bound~\eqref{second-assumption-in-homotopy-lemma} of Section~\ref{sec:radii-polynomial} plays the same role for $D\overline{F}(\tilde{w})$ in Remark~\ref{rem:border-u0}.

Lemma~\ref{lem:bordered} is a statement in the Hilbert space $\mathbb{C}^2 \times \ell^2$, whereas the computer-assisted proof gives information in $X_1 \times X_2$, for instance that $Q(\tilde{w})$ is injective there. For counting negative eigenvalues this difference does not matter: a subspace on which a quadratic form is negative can be approximated by finite Fourier sums, so it is also ``visible'' in $X_1 \times X_2$. For the nullspace, however, it does matter, since injectivity on $X_1 \times X_2$ does not exclude a kernel vector in $\ell^2$ outside $\ell^1_\nu$. The following lemma excludes such vectors. Since $DF(\tilde{u})$ is dominated by $L$ on the high modes, its eigenvectors inherit the decay of $\tilde{u}$.
\begin{lemma}\label{lem:regularity}
Let $\lambda \in \mathbb{R}$, and let $\phi \in \ell^2$ satisfy $DF(\tilde{u})\phi = \lambda\phi$. Then $\phi \in \ell^1_\nu$.
\end{lemma}
\begin{proof}
Write $DF(\tilde{u}) = L + V$, where $V$ is the convolution with $v \bydef -2\gamma\tilde{u} + 3\tilde{u}^2$. Since $\tilde{u} \in \ell^1_\nu$ and $\ell^1_\nu$ is a Banach algebra, $v \in \ell^1_\nu$; this is where the decay of $\tilde{u}$ enters. Hence $V$ is bounded with norm at most $c \bydef \|v\|_{\ell^1_\nu}$ on $\ell^1_\nu$ and on $\ell^2$, since $\|Vh\|_2 \leq \|v\|_1 \|h\|_2 \leq \|v\|_{\ell^1_\nu} \|h\|_2$ by Young's inequality. Choose $N'$ so large that $m_{N'} \bydef \min_{n \notin I^{N'}} |L_{nn} - \lambda| > c$. The truncation $\pi^{\le N'}\phi$ has only finitely many nonzero entries, hence it lies in $\ell^1_\nu$, so it remains to show that the tail $\pi^{>N'}\phi$ lies in $\ell^1_\nu$. To this end we characterize the tail as the fixed point of a contraction. On the modes outside $I^{N'}$ the eigenvalue equation reads $(L-\lambda)\pi^{>N'}\phi = -\pi^{>N'}V\phi$, so $\psi = \pi^{>N'}\phi$ is a fixed point of
\[
  \Theta(\psi) \bydef -(L-\lambda)^{-1}\pi^{>N'} V(\pi^{\le N'}\phi + \psi).
\]
Since $(L-\lambda)^{-1}$ acts diagonally outside $I^{N'}$, $\|(L-\lambda)^{-1}\pi^{>N'}\| \leq m_{N'}^{-1}$ on $\ell^2$ and on $\ell^1_\nu$, and the map $\Theta$ has Lipschitz constant at most $c/m_{N'} < 1$ on both $\pi^{>N'}\ell^2$ and $\pi^{>N'}\ell^1_\nu$. It maps both $\pi^{>N'}\ell^2$ and $\pi^{>N'}\ell^1_\nu$ into themselves, since $\pi^{\le N'}\phi \in \ell^1_\nu \subset \ell^2$. By the Banach fixed point theorem, $\Theta$ has a unique fixed point in $\pi^{>N'}\ell^2$, which is $\pi^{>N'}\phi$ since $\phi \in \ell^2$, and a fixed point $\tilde{\psi}$ in $\pi^{>N'}\ell^1_\nu$. Since $\ell^1_\nu \subset \ell^2$, we get that $\tilde{\psi}$ is also a fixed point in $\pi^{>N'}\ell^2$, and by uniqueness $\tilde{\psi} = \pi^{>N'}\phi$. Hence $\phi = \pi^{\le N'}\phi + \tilde{\psi} \in \ell^1_\nu$.
\end{proof}

We now combine the two lemmas. Applying Lemma~\ref{lem:bordered} to count negative eigenvalues of $DF(\tilde{u})$ needs two facts that are not known a priori: that $\partial_1\tilde{u}$ and $\partial_2\tilde{u}$ are linearly independent, and that $Q(\tilde{w})$ is injective on $\mathbb{C}^2 \times \ell^2$. The next lemma derives both from two conditions that the computer-assisted proof can verify. The first is that the matrix $M$ of Section~\ref{sec:beta-zero} is invertible. The second is the inequality~\eqref{condition-radii-polynomial-q} of Theorem~\ref{th:radii-polynomial-q}, which makes $\overline{A}Q(\tilde{w})$ invertible on $X_1 \times X_2$ and hence $Q(\tilde{w})$ injective there. Lemma~\ref{lem:regularity} then extends the injectivity from $X_1 \times X_2$ to $\mathbb{C}^2 \times \ell^2$.
\begin{lemma}\label{lem:eigenvalues-of-Q-2d}
Let the assumptions of Theorem~\ref{th:radii-polynomial-q} hold, and let the matrix $M$ in~\eqref{eq:def-C} be invertible. Then $\partial_1\tilde{u}$ and $\partial_2\tilde{u}$ are linearly independent, $Q(\tilde{w})$ is injective on $\mathbb{C}^2 \times \ell^2$, $\mathrm{Null}(DF(\tilde{u})) = K$,
\[
  \sigma(Q(\tilde{w})) = \bigl(\sigma(DF(\tilde{u})) \setminus \{0\}\bigr) \cup \{\pm\lambda_1,\pm\lambda_2\},
\]
and $\mathcal{N}(DF(\tilde{u})) = \mathcal{N}(Q(\tilde{w})) - 2$.
\end{lemma}
\begin{proof}
We first show that $\partial_1 \tilde{u}$ and $\partial_2 \tilde{u}$ are linearly independent. Since $\tilde{u}$ and $u_0$ are real-valued, $M$ is real. If $\alpha_1 \partial_1 \tilde{u} + \alpha_2 \partial_2 \tilde{u} = 0$, then $(M^T\alpha)_j = (\alpha_1 \partial_1 \tilde{u} + \alpha_2 \partial_2 \tilde{u}, \partial_j u_0)_2 = 0$, so $\alpha = 0$ since $M$ is invertible.

Next we prove that $Q(\tilde{w})$ is injective on $\mathbb{C}^2 \times \ell^2$. Namely, taking $(\beta,\phi) \in \mathbb{C}^2 \times \ell^2$ with $Q(\tilde{w})(\beta,\phi) = 0$ we set out to prove that $(\beta,\phi)=(0,0)$. Written out, $Q(\tilde{w})(\beta,\phi) = 0$ means
\begin{align}
\label{e:Q12}
  (\phi,\partial_1\tilde{u})_2 = (\phi,\partial_2\tilde{u})_2 = 0 \qquad \text{and} \qquad \beta_1\partial_1\tilde{u} + \beta_2\partial_2\tilde{u} + DF(\tilde{u})\phi = 0 .
\end{align}
We first show that $\beta = 0$. Take the inner product of the second equation with $\partial_i\tilde{u}$. The term with $DF(\tilde{u})\phi$ drops out, because $DF(\tilde{u})$ is self-adjoint and $DF(\tilde{u})\partial_i\tilde{u} = 0$, so that $(DF(\tilde{u})\phi,\partial_i\tilde{u})_2 = (\phi,DF(\tilde{u})\partial_i\tilde{u})_2 = 0$. What remains is $\sum_j \beta_j (\partial_j\tilde{u},\partial_i\tilde{u})_2 = (\Phi\beta)_i = 0$ for $i = 1,2$. Since $\partial_1\tilde{u}$ and $\partial_2\tilde{u}$ are linearly independent, $\Phi$ is invertible, so $\beta = 0$.

Now that we know $\beta=0$, the second equation in~\eqref{e:Q12} gives $DF(\tilde{u})\phi = 0$, so $\phi \in \ell^1_\nu$ by Lemma~\ref{lem:regularity}. Hence $(0,\phi)$ lies in $X_1 \times X_2$ and is mapped to zero by $Q(\tilde{w})$, where we have used the first equation in~\eqref{e:Q12}. By~\eqref{first-assumption-in-homotopy-lemma}, $\overline{A}Q(\tilde{w})$ is invertible on $X_1 \times X_2$, so $Q(\tilde{w})$ is injective there, and $\phi = 0$.

The remaining assertions follow directly from Lemma~\ref{lem:bordered}.
\end{proof}

Lemma~\ref{lem:eigenvalues-of-Q-2d} reduces the count for $DF(\tilde{u})$ to a count for $Q(\tilde{w})$, which is not known explicitly. By~\eqref{first-assumption-in-homotopy-lemma}, $\overline{A}$ is an approximate inverse of $Q(\tilde{w})$. The next proposition shows that, as a consequence, $\overline{A}$ and $Q(\tilde{w})$ have the same number of negative eigenvalues. Since $\overline{A} = \overline{A}^N \oplus \pi^{>N} L^{-1}$ has a positive tail, this is the number of negative eigenvalues of the matrix $\overline{A}^N$, which is determined rigorously in Section~\ref{sec:numerical-aspects}.
\begin{prop}\label{prop:homotopy}
Let the assumptions of Theorem~\ref{th:radii-polynomial-q} hold. Then $\mathcal{N}(Q(\tilde{w})) = \mathcal{N}(\overline{A}^N)$.
\end{prop}
\begin{proof}
This follows from a homotopy from $\overline{A}^{-1}$ to $Q(\tilde{w})$ along which~\eqref{first-assumption-in-homotopy-lemma} persists; see Appendix~\ref{apen:homotopy}.
\end{proof}
\begin{corollary}\label{lem:homotopy-arguments}
Let the assumptions of Theorem~\ref{th:radii-polynomial-q} hold and let $M$ be invertible. Then $\mathcal{N}(DF(\tilde{u})) = \mathcal{N}(\overline{A}^N) - 2$. In particular, $\tilde{u}$ is stable if and only if $\overline{A}^N$ has exactly two negative eigenvalues.
\end{corollary}
\begin{proof}
Combine Lemma~\ref{lem:eigenvalues-of-Q-2d} and Proposition~\ref{prop:homotopy}.
\end{proof}

The bound~\eqref{first-assumption-in-homotopy-lemma} gives more than the count of negative eigenvalues: it also yields a lower bound on the absolute values of the nonzero eigenvalues of $DF(\tilde{u})$. The idea is that a nonzero eigenvalue $\lambda$ of $DF(\tilde{u})$ is also an eigenvalue of $Q(\tilde{w})$ by Lemma~\ref{lem:eigenvalues-of-Q-2d}, with an eigenvector in $X_1 \times X_2$ by Lemma~\ref{lem:regularity}. Hence~\eqref{first-assumption-in-homotopy-lemma} and the Neumann series give $1/|\lambda| \leq \|Q(\tilde{w})^{-1}\|_\xi \leq \|\overline{A}\|_\xi/(1-\delta)$. This is the simplest bound, but the weighted column sums in $\|\overline{A}\|_\xi$ can be much larger than the largest absolute eigenvalue of $\overline{A}$, which is the natural size of $\|Q(\tilde{w})^{-1}\|$. We therefore estimate instead
\[
  1/|\lambda| \leq \|Q(\tilde{w})^{-1}\|_{\mathcal{B}(\ell^2)} \leq \|\overline{A}\|_{\mathcal{B}(\ell^2)} + \|Q(\tilde{w})^{-1} - \overline{A}\|_{\mathcal{B}(\ell^2)} .
\]
The error term is estimated through the weighted norm, by $\frac{\delta}{1-\delta}\|\overline{A}\|_\xi$, see the proof of Theorem~\ref{th:spectral-gap} below, so for small $\delta$ this improves the bound substantially. The following lemma, an instance of the Schur test, provides the transfer from $\|\cdot\|_\xi$ to $\ell^2$.
\begin{lemma}\label{lem:schur}
Let $E$ be a bounded operator on $X_1 \times X_2$ which is self-adjoint with respect to the inner product~\eqref{e:inner2}. Then $\|E\|_{\mathcal{B}(\ell^2)} \leq \|E\|_\xi$.
\end{lemma}
\begin{proof}
In the coordinates $(\beta_1,\beta_2,(a_n)_{n})$ the norm $\|\cdot\|_\xi$ is a weighted $\ell^1$ norm, with weights $\omega_k$ equal to $\xi_1$ for the components of $\beta$ and to $\xi_2 \nu_1^{|n_1|}\nu_2^{|n_2|}$ for the Fourier coefficients. The operator norm of a matrix $(E_{km})$ with respect to the weighted $\ell^1$ norm is the weighted maximal column sum $C \bydef \|E\|_\xi = \sup_m \sum_k |E_{km}| \omega_k/\omega_m$, so that
\begin{align}
\label{e:Ekm}
  \sum_k |E_{km}|\, \omega_k \leq C\,\omega_m \quad \text{for all } m, \qquad \sum_m |E_{km}|\, \omega_m \leq C\,\omega_k \quad \text{for all } k ,
\end{align}
where the second inequality follows from the first since $|E_{km}| = |E_{mk}|$ for self-adjoint $E$. By the Cauchy-Schwarz inequality and the second inequality in~\eqref{e:Ekm},
\[
  |(Eh)_k|^2 \leq \Bigl(\sum_m |E_{km}|\, \omega_m\Bigr) \Bigl(\sum_m |E_{km}| \frac{|h_m|^2}{\omega_m}\Bigr) \leq C\,\omega_k \sum_m |E_{km}| \frac{|h_m|^2}{\omega_m} .
\]
Summing over $k$ and using the first inequality in~\eqref{e:Ekm} gives $\|Eh\|_2^2 \leq C \sum_m |h_m|^2 \omega_m^{-1} \sum_k |E_{km}|\, \omega_k \leq C^2 \|h\|_2^2$.
\end{proof}
\begin{theorem}\label{th:spectral-gap}
Let the assumptions of Corollary~\ref{lem:homotopy-arguments} hold, and let $\delta < 1$ be as in~\eqref{first-assumption-in-homotopy-lemma}. Then every nonzero eigenvalue $\lambda$ of $DF(\tilde{u})$ satisfies
\[
    |\lambda| \geq \frac{1-\delta}{(1-\delta)\|\overline{A}\|_{\mathcal{B}(\ell^2)} + \delta\|\overline{A}\|_\xi} .
\]
\end{theorem}
\begin{proof}
By Lemma~\ref{lem:eigenvalues-of-Q-2d}, $Q(\tilde{w})$ is injective on $\mathbb{C}^2 \times \ell^2$; being self-adjoint with compact resolvent, it is invertible. We have $Q(\tilde{w})^{-1} = (\overline{A}Q(\tilde{w}))^{-1}\overline{A}$, and $\|(\overline{A}Q(\tilde{w}))^{-1}\|_\xi \leq 1/(1-\delta)$ by~\eqref{first-assumption-in-homotopy-lemma} and the Neumann series. Hence $\|Q(\tilde{w})^{-1}\|_\xi \leq \|\overline{A}\|_\xi/(1-\delta)$, and
\[
  \|Q(\tilde{w})^{-1} - \overline{A}\|_\xi = \|(I - \overline{A}Q(\tilde{w}))Q(\tilde{w})^{-1}\|_\xi \leq \frac{\delta}{1-\delta}\, \|\overline{A}\|_\xi .
\]
Since $Q(\tilde{w})^{-1} - \overline{A}$ is self-adjoint, Lemma~\ref{lem:schur} gives the same bound for $\|Q(\tilde{w})^{-1} - \overline{A}\|_{\mathcal{B}(\ell^2)}$. If $\lambda \neq 0$ is an eigenvalue of $DF(\tilde{u})$, then $\lambda$ is an eigenvalue of $Q(\tilde{w})$ by Lemma~\ref{lem:eigenvalues-of-Q-2d}, so $|\lambda|^{-1} \leq \|Q(\tilde{w})^{-1}\|_{\mathcal{B}(\ell^2)} \leq \|\overline{A}\|_{\mathcal{B}(\ell^2)} + \frac{\delta}{1-\delta} \|\overline{A}\|_\xi$.
\end{proof}
Both norms of $\overline{A}$ in Theorem~\ref{th:spectral-gap} are computable, since $\overline{A} = \overline{A}^N \oplus \pi^{>N} L^{-1}$. As $\overline{A}$ is self-adjoint, its norm in $\mathcal{B}(\ell^2)$ is the largest absolute value of its eigenvalues, and the weighted column sums in the proof of Lemma~\ref{lem:schur} give its norm in $\|\cdot\|_\xi$:
\[
  \|\overline{A}\|_{\mathcal{B}(\ell^2)} = \max\bigl(\rho(\overline{A}^N), L_N^{-1}\bigr), \qquad \|\overline{A}\|_\xi = \max\bigl(\|\overline{A}^N\|_\xi, L_N^{-1}\bigr).
\]
Here the spectral radius $\rho(\overline{A}^N)$ is the largest absolute value of the eigenvalues of the matrix $\overline{A}^N$, $\|\overline{A}^N\|_\xi$ is its largest weighted column sum, and $L_N^{-1}$ is the largest entry of the tail $\pi^{>N} L^{-1}$, with $L_N$ as in~\eqref{def:L_N}.

For a stable solution, Theorem~\ref{th:spectral-gap} makes the stability statement of Section~\ref{sec:variational} quantitative. All nonzero eigenvalues of $DF(\tilde{u})$ are then positive and at least the bound in the theorem, so under the linearized flow $\partial_t v = -DF(\tilde{u})v$ perturbations orthogonal to $\partial_1\tilde{u}$ and $\partial_2\tilde{u}$ decay at least at this exponential rate. For an unstable solution, the perturbations in the directions of negative eigenvalues grow at least at this rate.

\begin{remark}\label{rem:border-eigenvalues}
The bound of Theorem~\ref{th:spectral-gap} applies to all eigenvalues of $Q(\tilde{w})$, including the eigenvalues $\pm\lambda_1,\pm\lambda_2$ of Lemma~\ref{lem:bordered}, which come from the border and not from $DF(\tilde{u})$. Since $\overline{A}$ approximates $Q(\tilde{w})^{-1}$, the norm $\|\overline{A}\|_{\mathcal{B}(\ell^2)}$ is at least about $1/\min(\lambda_1,\lambda_2)$, so the bound cannot exceed about $\min(\lambda_1,\lambda_2)$, whatever the actual gap is. There are two ways to reduce this effect. First, the scale of the border is arbitrary. Replacing $\partial_j u$ by $s\,\partial_j u$ in $Q$, with $s > 0$, multiplies $\pm\lambda_1,\pm\lambda_2$ by $s$ and leaves the other eigenvalues unchanged. With the approximate inverse $D^{-1}\overline{A}D^{-1}$, where $D = \mathrm{diag}(s,s,I)$, and the weights $(s\xi_1,\xi_2)$, the bound $\delta$ is unchanged, since only the ratio $\xi_1/\xi_2$ enters. Taking $s$ large removes the eigenvalues near $\pm 1/\lambda_1,\pm 1/\lambda_2$ from the term with $\|\cdot\|_{\mathcal{B}(\ell^2)}$ in the bound, but not their contribution to the term with $\|\cdot\|_\xi$. Second, the eigenvalues of $Q(\tilde{w})^{-1}$ and $\overline{A}$ can be matched one by one, which makes it possible to exclude the four border eigenvalues explicitly, see Section~\ref{sec:conclusion}. It is through $Q$ that such statements about the eigenvalues become possible. The nullspace and the number of negative eigenvalues can also be obtained without $Q$, see Remark~\ref{rem:border-u0}.
\end{remark}

\begin{remark}\label{rem:border-u0}
The conclusions $\mathrm{Null}(DF(\tilde{u})) = \mathrm{span}\{\partial_1\tilde{u},\partial_2\tilde{u}\}$ of Lemma~\ref{lem:eigenvalues-of-Q-2d} and $\mathcal{N}(DF(\tilde{u})) = \mathcal{N}(\overline{A}^N) - 2$ of Corollary~\ref{lem:homotopy-arguments} can also be obtained from $D\overline{F}(\tilde{w})$ itself, without the bound $Z_Q$. The counting argument is less direct than for $Q$, since the border $\partial_j u_0$ of $D\overline{F}(\tilde{w})$ does not consist of null vectors of $DF(\tilde{u})$, so the decoupling in the proof of Lemma~\ref{lem:bordered} is not available. Instead, a change of coordinates makes the quadratic form block diagonal, and for this only the invertibility of $M$ is needed.

Let $K_0 \bydef \mathrm{span}\{\partial_1 u_0, \partial_2 u_0\}$. Every $x \in \ell^2$ can be written uniquely as $x = \alpha_1\partial_1\tilde{u} + \alpha_2\partial_2\tilde{u} + y$ with $\alpha \in \mathbb{C}^2$ and $y \in K_0^\perp$, since $(x,\partial_i u_0)_2 = (M^T\alpha)_i$ determines $\alpha$ when $M$ is invertible. Hence
\[
  J(\alpha,y) \bydef \alpha_1\partial_1\tilde{u} + \alpha_2\partial_2\tilde{u} + y, \qquad \overline{J}(\beta,\alpha,y) \bydef \bigl(\beta, J(\alpha,y)\bigr)
\]
are bijections from $\mathbb{C}^2 \times K_0^\perp$ onto $\ell^2$ and from $\mathbb{C}^2 \times \mathbb{C}^2 \times K_0^\perp$ onto $\mathbb{C}^2 \times \ell^2$. Since $DF(\tilde{u})\partial_j\tilde{u} = 0$, $(\partial_i\tilde{u},\partial_j u_0)_2 = M_{ij}$ and $(y,\partial_i u_0)_2 = 0$, a direct computation gives
\[
  J^* DF(\tilde{u})\, J = \begin{bmatrix} 0 & 0 \\ 0 & D^\perp \end{bmatrix}, \qquad \overline{J}^* D\overline{F}(\tilde{w})\, \overline{J} = \begin{bmatrix} 0 & M^T & 0 \\ M & 0 & 0 \\ 0 & 0 & D^\perp \end{bmatrix},
\]
where $D^\perp$ is $DF(\tilde{u})$ restricted to $K_0^\perp$ in the sense of quadratic forms: $(D^\perp y,z)_2 = (DF(\tilde{u})y,z)_2$ for $y,z \in K_0^\perp$.

The maps $J$ and $\overline{J}$ are not unitary, so the block forms do not have the same eigenvalues as $DF(\tilde{u})$ and $D\overline{F}(\tilde{w})$. They do have the same number of negative eigenvalues, by Sylvester's law of inertia, and this number is all we transfer. By the min-max principle, $\mathcal{N}(DF(\tilde{u}))$ is the largest dimension of a subspace on which the quadratic form $(DF(\tilde{u})x,x)_2$ is negative for all $x \neq 0$. Since $(J^*DF(\tilde{u})Jz,z)_2 = (DF(\tilde{u})Jz,Jz)_2$, the bijection $J$ maps such subspaces for the block form onto such subspaces for $DF(\tilde{u})$, and similarly for $\overline{J}$. For the block forms and for $D^\perp$ we define $\mathcal{N}$ in the same way.

For a block diagonal form, $\mathcal{N}$ is the sum over the blocks. The $4 \times 4$ block with off-diagonal blocks $M^T$ and $M$ has the eigenvalues $\pm\sigma_1,\pm\sigma_2$, where $\sigma_1^2,\sigma_2^2 > 0$ are the eigenvalues of $M^TM$, so two of its eigenvalues are negative. Hence $\mathcal{N}(DF(\tilde{u})) = \mathcal{N}(D^\perp)$ and $\mathcal{N}(D\overline{F}(\tilde{w})) = \mathcal{N}(D^\perp) + 2$, so that $\mathcal{N}(D\overline{F}(\tilde{w})) = \mathcal{N}(DF(\tilde{u})) + 2$.

It remains to relate $\mathcal{N}(D\overline{F}(\tilde{w}))$ to $\mathcal{N}(\overline{A}^N)$ and to determine the nullspace of $DF(\tilde{u})$. Both use the bound~\eqref{second-assumption-in-homotopy-lemma} of the existence proof, in the norm $\|\cdot\|_\chi$. The argument of Appendix~\ref{apen:homotopy} applies with $D\overline{F}(\tilde{w})$ in place of $Q(\tilde{w})$ and~\eqref{second-assumption-in-homotopy-lemma} in place of~\eqref{first-assumption-in-homotopy-lemma}, see Remark~\ref{rem:homotopy-DFbar}, so $\mathcal{N}(D\overline{F}(\tilde{w})) =\mathcal{N}(\overline{A}^N)$ and hence $\mathcal{N}(DF(\tilde{u})) = \mathcal{N}(\overline{A}^N) - 2$.

For the nullspace, let $DF(\tilde{u})x = 0$ and write $x = \alpha_1\partial_1\tilde{u} + \alpha_2\partial_2\tilde{u} + y$ with $y \in K_0^\perp$. Then $DF(\tilde{u})y = 0$, so $y \in \ell^1_\nu$ by Lemma~\ref{lem:regularity}, and $D\overline{F}(\tilde{w})(0,y) = 0$. Since $D\overline{F}(\tilde{w})$ is injective on $X_1 \times X_2$ by~\eqref{second-assumption-in-homotopy-lemma}, $y = 0$. Hence $\mathrm{Null}(DF(\tilde{u})) = \mathrm{span}\{\partial_1\tilde{u},\partial_2\tilde{u}\}$, which is two-dimensional since $M$ is invertible.
\end{remark}

\section{Rigorous numerics}
\label{sec:bounds-2D}

We construct the approximate solution $u_0$ and the approximate inverse $\overline{A}$, explain how $\mathcal{N}(\overline{A})$ is determined, and compute the bounds of Theorems~\ref{th:radii-polynomial} and~\ref{th:radii-polynomial-q}.

\subsection{Numerical Aspects}
\label{sec:numerical-aspects}
Recall the index set $I^N$ and the projections $\pi^{\le N}$ and $\pi^{>N}$ from Section~\ref{sec:radii-polynomial}. We also define
\begin{align}
    \overline{\pi}^{\le N} \bydef \begin{bmatrix}
        I & 0 \\
        0 & \pi^{\le N}
    \end{bmatrix},~~\overline{\pi}^{>N} \bydef \begin{bmatrix}
        0 & 0 \\
        0 & \pi^{>N}
    \end{bmatrix}\label{def:projection-operators}
\end{align}
which will be useful when working on the product space $X_1 \times X_2$. To obtain a numerical approximation of the solution, we first fix the domain $[-d_1,d_1] \times [-d_2,d_2]$. We choose $d_1$ and $d_2$ small, since without symmetry we work with full Fourier series with complex coefficients, and the cost of storing and manipulating the corresponding matrices grows quickly with the domain size. We then construct an initial guess as a sum of functions of the form $c \mathrm{sech}(\sqrt{(x_1 + a)^2 + (x_2 + b)^2})^2$ for some parameters $c > 0$ and any $a,b$. For example,
\begin{align}
    \overline{u}_{1}(x) = \sum_{i = 1}^l c_i \mathrm{sech}(\sqrt{(x_1 + a_i)^2 + (x_2 + b_i)^2})^2
\end{align}
for some $c_i,a_i,b_i$ and number of spots $l$. By changing the values of the parameters, specifically $a_i$ and $b_i$, we can place $l$ "spots" in different parts of the domain. We randomly placed these spots, and then compute a Fourier series approximation of the function $\overline{u}_1$. That is, we compute $(\overline{a}_n)_{n \in \mathbb{Z}^2}$ such that
\begin{align}
    \overline{u}_1(x) \approx \sum_{n \in \mathbb{Z}^2} \overline{a}_n \mathrm{e}^{i \tilde{n} \cdot x} \bydef \overline{u}_0(x).
\end{align}
From here, we run a Newton method on $\overline{u}_0$ to refine our initial guess. If Newton converges, then we denote this solution $u_0$ and use it as our approximate solution for the proof. Our approach is not guaranteed to work for any PDE and solution; however, it was sufficient for our purposes as Newton converged. 
\par The approximate inverse $\overline{A}$ consists of two parts: a finite part which can be represented as a matrix, and an infinite part which we control. As $\overline{A}$ is supposed to be an approximate inverse for $D\overline{F}(w_0)$, we choose
\begin{align}
    \overline{A}^N \approx (\overline{\pi}^{\le N} D\overline{F}(w_0)\overline{\pi}^{\le N})^{-1}.
\end{align}
The operator $\overline{A}^N$ can be represented as a matrix, which can be stored on the computer. This is the finite part of $\overline{A}$. For the tail, rather than considering the full $D\overline{F}(w_0)^{-1}$, we only choose the inverse of the linear part, $L$. Note that $L$ may not be invertible, so we abuse notation and write
\begin{align}
    (\pi^{>N} L^{-1})_{nm} \bydef \begin{cases}
        \frac{\delta_{nm}}{L_{nn}} & n,m \in \mathbb{Z}^2 \setminus I^N \\
        0 & \mathrm{else}
    \end{cases},
\end{align}
where $\delta_{nm}$ is the Kronecker delta. Then, we define
\begin{align}
    \overline{A} \bydef \overline{A}^N + \begin{bmatrix}
        0 & 0 \\
        0 & \pi^{>N} L^{-1}
    \end{bmatrix} = \begin{bmatrix}
        \overline{A}_{11} & \overline{A}_{12} \\
        \overline{A}_{21} & \overline{A}_{22}
    \end{bmatrix},~~
    \overline{\pi}^{>N} \overline{A} = \begin{bmatrix}
       0 & 0 \\
        0 & \pi^{>N} L^{-1}
    \end{bmatrix}
\end{align}
where $\overline{A}_{ij} \bydef P_i \overline{A} P_j \in \mathcal{B}(Y_j,X_i)$. We replace the numerically computed matrix $\overline{A}^N$ by $\frac{1}{2}(\overline{A}^N + (\overline{A}^N)^*)$, so that $\overline{A}$ is self-adjoint with respect to the inner product~\eqref{e:inner2}, as required in Section~\ref{sec:stability-2D}. Also, we will choose
\begin{align}
    \overline{A}_{12} = \overline{A}_{12} \pi^{\le N}, ~ \overline{A}_{21} = \pi^{\le N} \overline{A}_{21}.
\end{align}
Based on how we defined $\overline{A}$ earlier, it is also natural to define
\begin{align}
    \overline{A}_{22}^N \bydef \pi^{\le N} \overline{A}_{22},~~ \pi^{>N} \overline{A}_{22} = \pi^{>N} L^{-1},~~ \overline{A}^N = \begin{bmatrix}
        \overline{A}_{11} & \overline{A}_{12} \\
        \overline{A}_{21} & \overline{A}_{22}^N
    \end{bmatrix} = \begin{bmatrix}
        \overline{A}_{11} & \overline{A}_{12}\pi^{\le N} \\
        \pi^{\le N} \overline{A}_{21} & \overline{A}_{22}^N
    \end{bmatrix}.
\end{align}
Our choices will also allow us to bound the tail of $\overline{A}$. In fact,
\begin{align}
    \|\overline{\pi}^{>N} \overline{A}\|_{\mathcal{B}(X_1 \times X_2)} = \|\pi^{>N} \overline{A}_{22}\|_{\mathcal{B}(X_2)} &= \|\pi^{>N} L^{-1}\|_{\mathcal{B}(X_2)} \leq \frac{1}{L_N},
\end{align}
where
\begin{align}
    L_N \bydef \min_{n \in \mathbb{Z}^2 \setminus I^N} |L_{nn}| = \bigl(\zeta_N^2 - 1\bigr)^2 + \mu \qquad \text{if } \zeta_N \geq 1, \qquad \zeta_N \bydef \min_{i=1,2} \frac{(N_i+1)\pi}{d_i} .\label{def:L_N}
\end{align}
Here $\zeta_N$ is the smallest value of $|\tilde{n}|$ for $n \notin I^N$, and we used that $L_{nn} > 0$ for $n \notin I^N$ by the choice of $N$ in Section~\ref{sec:radii-polynomial}.
This bound is used frequently below.
\par Since $\overline{A}_{12} = \overline{A}_{12}\pi^{\le N}$ and $\overline{A}_{21} = \pi^{\le N}\overline{A}_{21}$, the operator $\overline{A}$ is block diagonal, $\overline{A} = \overline{A}^N \oplus \pi^{>N} L^{-1}$. For Corollary~\ref{lem:homotopy-arguments} and Theorem~\ref{th:spectral-gap} it remains to determine $\mathcal{N}(\overline{A}^N)$ and to bound $\rho(\overline{A}^N)$. To this end we use Gershgorin's theorem, as in~\cite{JB_symmetries_2}. Let $V$ be a matrix whose columns are numerically computed approximate eigenvectors of $\overline{A}^N$, and let
\[
    \overline{\mathcal{A}}^N \bydef V^{-1} \overline{A}^N V .
\]
We compute $V^{-1}$ rigorously by inverting the interval matrix $V$. We refer to $\overline{\mathcal{A}}^N$ as the pseudo-diagonalization of $\overline{A}^N$. It is nearly diagonal, and it has the same eigenvalues as $\overline{A}^N$, which are real since $\overline{A}^N$ is self-adjoint. Write $\overline{\mathcal{A}}^N = (a_{ik})$, and let $D_i$ be the closed disc in $\mathbb{C}$ with center $a_{ii}$ and radius $\sum_{k \neq i} |a_{ik}|$. If $|\mathrm{Re}\, a_{ii}| > \sum_{k \neq i} |a_{ik}|$ for all $i$, so that no disc meets the imaginary axis, then by Gershgorin's theorem~\cite[Section~6.1]{horn_johnson_2013} $\mathcal{N}(\overline{A}^N)$ equals the number of discs in the left half-plane. Moreover, $\rho(\overline{A}^N) \leq \max_i \sum_k |a_{ik}| = \|\overline{\mathcal{A}}^N\|_\infty$, where $\|\cdot\|_\infty$ is the maximal absolute row sum.

\subsection{Computing the Bounds}

We compute the bounds of Theorem~\ref{th:radii-polynomial}, beginning with $Y^{[1]}$ and $Y^{[2]}$.
\begin{lemma}\label{lem:Y0-2d}
Let $Y^{[1]}, Y^{[2]} \geq 0$ be defined as
\begin{align}
    &Y^{[1]} \bydef \|P_1 \overline{A}^N \overline{F}(w_0)\|_{X_1} \\
    &Y^{[2]} \bydef \|P_2 \overline{A}^N \overline{F}(w_0)\|_{X_2} + \frac{1}{L_N}\|(\pi^{\le N_0} - \pi^{\le N}) Lu_0 + (\pi^{\le 3N_0} - \pi^{\le N}) G(u_0)\|_{\ell^1_\nu}
\end{align}
where $L_N$ is defined in \eqref{def:L_N}.
Then $\|P_i [T(w_0) - w_0]\|_{X_i} \leq Y^{[i]}$ for $i = 1,2$.
\end{lemma} 
\begin{proof}
The bound $Y^{[1]}$ can be directly evaluated on the computer. For $Y^{[2]}$, we start from the definition.
\begin{align}
    \|P_2[T(w_0) - w_0]\|_{X_2} &= \|P_2 \overline{A}\,\overline{F}(w_0)\|_{X_2} \\
    &= \|\pi^{\le N} P_2 \overline{A}\,\overline{F}(w_0)\|_{X_2} + \|\pi^{>N} P_2 \overline{A}\,\overline{F}(w_0)\|_{X_2}\\
    &= \|P_2 \overline{A}^N \overline{F}(w_0)\|_{X_2} + \|\pi^{>N} (\overline{A}_{21} P_1 \overline{F}(w_0) + \overline{A}_{22} P_2 \overline{F}(w_0))\|_{X_2} \end{align}
The first term above can be evaluated on the computer. For the second term, recall that $\pi^{>N} \overline{A}_{21} = 0$ and $\|\pi^{>N} \overline{A}_{22}\|_{\mathcal{B}(X_2)} \leq \frac{1}{L_N}$. Then, 
\begin{align}
    \|\pi^{>N} (\overline{A}_{21} P_1 \overline{F}(w_0) + \overline{A}_{22} P_2 \overline{F}(w_0))\|_{X_2}  &=\|\pi^{>N} \overline{A}_{22}P_2 \overline{F}(w_0)\|_{X_2} \\
    &\leq \|\pi^{>N} \overline{A}_{22}\|_{\mathcal{B}(X_2)} \|\pi^{>N} P_2 \overline{F}(w_0)\|_{X_2}\\
    &\leq \frac{1}{L_N} \|\pi^{>N} F(u_0)\|_{\ell^1_\nu} \\
    &= \frac{1}{L_N}\|(\pi^{\le N_0} - \pi^{\le N}) L u_0 + (\pi^{\le 3N_0} - \pi^{\le N}) G(u_0)\|_{\ell^1_\nu}
\end{align}
where we used the definition of $L$ and $G$ from \eqref{def:L-and-G-sh-2d}.
As a result, we see that
\begin{align}
    \|P_2[T(w_0) - w_0]\|_{X_2} &\leq \|P_2 \overline{A}^N \overline{F}(w_0)\|_{X_2} + \frac{1}{L_N}\|(\pi^{\le N_0} - \pi^{\le N}) L u_0 + (\pi^{\le 3N_0} - \pi^{\le N}) G(u_0)\|_{\ell^1_\nu}
\end{align}
as desired.
\end{proof}
For a sequence $a \in \ell^1_\nu$ we denote by $\mathbb{a}$ the multiplication operator by the function with Fourier coefficients $a$, that is, $\mathbb{a}h \bydef a * h$. For instance, $DG(u_0) = -2\gamma\mathbb{u}_0 + 3\mathbb{u}_0^2 = \mathbb{v}_0$ with $v_0 \bydef -2\gamma u_0 + 3u_0^2$. The bounds $W_{22}^{[i]}$ are as follows; as explained after Theorem~\ref{th:radii-polynomial}, the other $W_{jk}^{[i]}$ vanish. We write $D_u\overline{F}(w) \bydef D\overline{F}(w)P_2$ for the derivative with respect to $u$.
\begin{lemma}\label{lem:Z2}
Let $q = (q_n)_{n \in \mathbb{Z}^2} = (-2\gamma \delta_n + 6a_n)_{n \in \mathbb{Z}^2}$ where $u_0 = (a_n)_{n \in \mathbb{Z}^2}$ and $\delta_n$ is the Kronecker delta. Let $L_N$ be defined as in \eqref{def:L_N}, and let
\begin{align*}
    &W_{22}^{[1]} \bydef \left\|\overline{A}_{12}\mathbb{q}\pi^{\le N+N_0}
    \right\|_{\mathcal{B}(X_2,X_1)} + 3\|\overline{A}_{12}\|_{\mathcal{B}(X_2,X_1)}r_2^* , \\
    &W_{22}^{[2]} \bydef \left\|\overline{A}_{22}^N\mathbb{q}\pi^{\le N+N_0}
    \right\|_{\mathcal{B}(X_2)}+\frac{1}{L_N} \|q\|_{\ell^1_\nu} + 3\left(\|\overline{A}_{22}^N\|_{\mathcal{B}(X_2)} + \frac{1}{L_N}\right)r_2^* .
\end{align*}
Then for all $0 < r \leq r^*$, $w \in B_r(w_0)$ and $i = 1,2$
\[
    \|P_i\overline{A}(D_u\overline{F}(w) - D_u\overline{F}(w_0))\|_{\mathcal{B}(X_2,X_i)} \leq W_{22}^{[i]} r_2 .
\]
\end{lemma}
\begin{proof}
The proof can be found in Appendix \ref{apen:Z2}.
\end{proof}
For the bounds $Z_j^{[i]}$ let $u_0^N \bydef \pi^{\le N} u_0$ and $v_0^N \bydef \pi^{\le N} v_0$. Note that $u_0$ and $v_0$ have finitely many nonzero Fourier coefficients, with indices in $I^{N_0}$ and $I^{2N_0}$, respectively. By performing these truncations, we are able to compute a sharp estimate with only $N$ coefficients that is less computationally expensive. The remaining terms are still finite, and we can perform a rougher estimate. This rough estimate is ok as long as the decay in $u_0$ and $v_0$ is good up to $N$. We also write $\partial u_0 \bydef (\partial_1 u_0, \partial_2 u_0)$, viewed as the map $\gamma \mapsto \gamma_1 \partial_1 u_0 + \gamma_2 \partial_2 u_0$ from $X_1$ to $X_2$, with $(\partial u_0)^* v = ((v,\partial_1 u_0)_2,(v,\partial_2 u_0)_2)$, and similarly $\partial u_0^N$. Then $D\overline{F}(w_0) = \begin{bmatrix} 0 & 0 \\ 0 & L \end{bmatrix} + \overline{\mathbb{v}}_0$, where
\[
    \overline{\mathbb{v}}_0 \bydef \begin{bmatrix}
        0 & (\partial u_0)^* \\
        \partial u_0 & \mathbb{v}_0
    \end{bmatrix}, \qquad
    \overline{\mathbb{v}}_0^N \bydef \begin{bmatrix}
        0 & (\partial u_0^N)^* \\
        \partial u_0^N & \mathbb{v}_0^N
    \end{bmatrix}.
\]
These definitions will allow us to split operators into those of size $N$ and those of size $N_0$.

\begin{lemma}\label{lem:split-Z_1}
Let the following be non-negative constants 
\begin{align}
    &Z_1^{[1],N} \bydef \|I - \overline{A}_{12} \partial u_0^N\|_{\mathcal{B}(X_1)},\\
    &Z_2^{[1],N} \bydef \left\|\overline{A}_{11} (\partial u_0^N)^* + \overline{A}_{12} (L + \mathbb{v}_0^N)\right\|_{\mathcal{B}(X_2,X_1)}, \\
    &Z_1^{[2],N} \bydef \|\overline{A}_{22} \partial u_0^N\|_{\mathcal{B}(X_1,X_2)} \\
    &Z_2^{[2],N} \bydef \|I - \overline{A}_{21} (\partial u_0^N)^* - \overline{A}_{22} (L + \mathbb{v}_0^N)\|_{\mathcal{B}(X_2)}\\
    &Z_{1}^{[1],\infty} \bydef \|\overline{A}_{12}\|_{\mathcal{B}(X_2,X_1)}(\|\partial_1 (u_0 - u_0^N)\|_{\ell^1_\nu} + \|\partial_2 (u_0 - u_0^N)\|_{\ell^1_\nu})  \\
    &Z_{2}^{[1],\infty} \bydef \|\overline{A}_{11}\|_{\mathcal{B}(X_1)} (\|\partial_1 (u_0 - u_0^N)\|_{2} + \|\partial_2 (u_0 - u_0^N)\|_{2}) + \|\overline{A}_{12}\|_{\mathcal{B}(X_2,X_1)}\|v_0 - v_0^N\|_{\ell^1_\nu} \\
    &Z_{1}^{[2],\infty} \bydef \left(\|\overline{A}_{22}^N\|_{\mathcal{B}(X_2)} + \frac{1}{L_N}\right)(\|\partial_1 (u_0 - u_0^N)\|_{\ell^1_\nu} + \|\partial_2 (u_0 - u_0^N)\|_{\ell^1_\nu})\\
    &Z_{2}^{[2],\infty} \bydef \|\overline{A}_{21}\|_{\mathcal{B}(X_1,X_2)} (\|\partial_1 (u_0 - u_0^N)\|_{2} + \|\partial_2 (u_0 - u_0^N)\|_{2}) + \left(\|\overline{A}_{22}^N\|_{\mathcal{B}(X_2)} + \frac{1}{L_N}\right)\|v_0 - v_0^N\|_{\ell^1_\nu}
\end{align}
and let $Z_{j}^{[i]}$ for $i,j = 1,2$ be defined as
\begin{align}
    &Z_j^{[i]} \bydef Z_j^{[i],N} + Z_j^{[i],\infty}.
\end{align}
Then, it follows that the $Z_j^{[i]}$'s satisfy their respective estimates in the statement of Theorem \ref{th:radii-polynomial}.
\end{lemma}
\begin{proof}
We have
\[
    DT(w_0) = I - \overline{A}\left(\begin{bmatrix}
        0 & 0 \\
        0 & L
    \end{bmatrix} + \overline{\mathbb{v}}_0^N\right) - \overline{A}(\overline{\mathbb{v}}_0 - \overline{\mathbb{v}}_0^N).
\]
Multiplying out the blocks, the first operator equals
\[
    \begin{bmatrix}
        I - \overline{A}_{12}\, \partial u_0^N & -\overline{A}_{11} (\partial u_0^N)^* - \overline{A}_{12} (L + \mathbb{v}_0^N) \\
        -\overline{A}_{22}\, \partial u_0^N & I - \overline{A}_{21} (\partial u_0^N)^* - \overline{A}_{22} (L + \mathbb{v}_0^N)
    \end{bmatrix},
\]
whose blocks have norms $Z_j^{[i],N}$. With $u_0^\infty \bydef u_0 - u_0^N$, the second operator equals
\[
    \begin{bmatrix}
        \overline{A}_{12}\, \partial u_0^\infty & \overline{A}_{11} (\partial u_0^\infty)^* + \overline{A}_{12} (\mathbb{v}_0 - \mathbb{v}_0^N) \\
        \overline{A}_{22}\, \partial u_0^\infty & \overline{A}_{21} (\partial u_0^\infty)^* + \overline{A}_{22} (\mathbb{v}_0 - \mathbb{v}_0^N)
    \end{bmatrix}.
\]
Its blocks are bounded by $Z_j^{[i],\infty}$, using $\|\partial u_0^\infty\,\gamma\|_{X_2} \leq (\|\partial_1 u_0^\infty\|_{\ell^1_\nu} + \|\partial_2 u_0^\infty\|_{\ell^1_\nu})\|\gamma\|_{X_1}$, $\|(\partial u_0^\infty)^* v\|_{X_1} \leq (\|\partial_1 u_0^\infty\|_2 + \|\partial_2 u_0^\infty\|_2)\|v\|_{X_2}$, $\|\mathbb{v}_0 - \mathbb{v}_0^N\|_{\mathcal{B}(X_2)} \leq \|v_0 - v_0^N\|_{\ell^1_\nu}$ and $\|\overline{A}_{22}\|_{\mathcal{B}(X_2)} \leq \|\overline{A}_{22}^N\|_{\mathcal{B}(X_2)} + \frac{1}{L_N}$.
\end{proof}

\begin{remark}
Lemma \ref{lem:split-Z_1} allows us to use $N_0$ for the size of sequences and $N$ for the size of operators. This is particularly useful in 2D, where one of the numerical challenges can be memory. By using different sizes for operators and sequences, we can gain the benefit of having a larger truncation, but also smaller operators to avoid memory problems.   
\end{remark}
The bounds $Z_j^{[i],N}$ are computed as follows.
\begin{lemma}\label{lem:Z1N-computation}
Let $L_N$ be defined as in \eqref{def:L_N}. Let $Z_j^{[i],N} \geq 0$ be defined as
\begin{align}
    &Z_1^{[1],N} \bydef \|I - \overline{A}_{12} \partial u_0^N\|_{\mathcal{B}(X_1)},~~Z_{2}^{[1],N} \bydef \|(\overline{A}_{11}(\partial u_0^N)^* + \overline{A}_{12} (L + \mathbb{v}_0^N)) \pi^{\le 2N}\|_{\mathcal{B}(X_2,X_1)} \\
    &Z_1^{[2],N} \bydef \| \overline{A}_{22}^N \partial u_0^N\|_{\mathcal{B}(X_1,X_2)},~~Z_{2}^{[2],N} \bydef \|\pi^{\le N}( I - \overline{A}_{21} (\partial u_0^N)^* - \overline{A}_{22} (L + \mathbb{v}_0^N)) \pi^{\le 2N}\|_{\mathcal{B}(X_2)} + \frac{1}{L_N} \|v_0^N\|_{\ell^1_\nu}.
\end{align}
Then, the $Z_{j}^{[i],N}$'s satisfy the bounds stated in Lemma \ref{lem:split-Z_1}.
\end{lemma}
\begin{proof}
The proof can be found in Appendix \ref{apen:Z1N}.
\end{proof}

\subsection{Computation of \texorpdfstring{$Z_{Q,j}^{[i]}$}{ZQj\^i}}
\label{sec:Zq}

The bounds~\eqref{eq:definition-Zq} in Theorem~\ref{th:radii-polynomial-q} measure the difference between the borders of $Q(w)$ and $D\overline{F}(w)$. For $w = (\beta,u) \in B_r(w_0)$ let $z \bydef u - u_0$, so that $\|z\|_{\ell^1_\nu} \leq r_2$, and let $\partial z \bydef (\partial_1 z, \partial_2 z)$, viewed as the map $\gamma \mapsto \gamma_1 \partial_1 z + \gamma_2 \partial_2 z$ from $X_1$ to $Y_2$, with $(\partial z)^* v = ((v,\partial_1 z)_2,(v,\partial_2 z)_2)$. Then
\[
  Q(w) - D\overline{F}(w) = \begin{bmatrix}
    0 & (\partial z)^* \\
    \partial z & 0
  \end{bmatrix},
  \qquad
  \overline{A}\bigl(Q(w) - D\overline{F}(w)\bigr) = \begin{bmatrix}
    \overline{A}_{12}\, \partial z & \overline{A}_{11} (\partial z)^* \\
    \overline{A}_{22}\, \partial z & \overline{A}_{21} (\partial z)^*
  \end{bmatrix}.
\]
\begin{lemma}\label{lem:Zq}
Let $\mathcal{D}$ be the diagonal operator with $\mathcal{D}_{nn} \bydef |\tilde{n}_1| + |\tilde{n}_2|$ for $n \neq 0$ and $\mathcal{D}_{00} \bydef 1$. Assume that $\zeta_N^2 \geq 1 + |\mu|$, and let $\kappa$ be as in Lemma~\ref{lem:derivative-estimate}. Define
\begin{align*}
    &Z_{Q,1}^{[1]}(r) \bydef \|\overline{A}_{12} \mathcal{D}\|_{\mathcal{B}(X_2,X_1)}\, r_2, &&Z_{Q,2}^{[1]}(r) \bydef 2\kappa\, \|\overline{A}_{11}\|_{\mathcal{B}(X_1)}\, r_2, \\
    &Z_{Q,1}^{[2]}(r) \bydef \Bigl(\|\overline{A}_{22}^N \mathcal{D}\|_{\mathcal{B}(X_2)} + \frac{\sqrt{2}\,\zeta_N}{L_N}\Bigr) r_2, &&Z_{Q,2}^{[2]}(r) \bydef 2\kappa\, \|\overline{A}_{21}\|_{\mathcal{B}(X_1,X_2)}\, r_2 .
\end{align*}
Then~\eqref{eq:definition-Zq} holds.
\end{lemma}
\begin{proof}
For $\gamma \in X_1$, the coefficients of $\mathcal{D}^{-1}\partial z\, \gamma$ are $i(\tilde{n}_1\gamma_1 + \tilde{n}_2\gamma_2) z_n / (|\tilde{n}_1| + |\tilde{n}_2|)$ for $n \neq 0$ and $0$ for $n = 0$, so $\|\mathcal{D}^{-1}\partial z\, \gamma\|_{X_2} \leq \|\gamma\|_{X_1} \|z\|_{\ell^1_\nu} \leq \|\gamma\|_{X_1}\, r_2$. Writing $\overline{A}_{j2}\, \partial z = (\overline{A}_{j2}\mathcal{D})(\mathcal{D}^{-1}\partial z)$ gives the bounds in the first column, once we note that
\[
  \|\overline{A}_{22}\mathcal{D}\|_{\mathcal{B}(X_2)} \leq \|\overline{A}_{22}^N\mathcal{D}\|_{\mathcal{B}(X_2)} + \sup_{n \notin I^N} \frac{\mathcal{D}_{nn}}{L_{nn}} \leq \|\overline{A}_{22}^N\mathcal{D}\|_{\mathcal{B}(X_2)} + \frac{\sqrt{2}\,\zeta_N}{L_N}.
\]
Here we used $\mathcal{D}_{nn} \leq \sqrt{2}\,|\tilde{n}|$, $|\tilde{n}| \geq \zeta_N$ for $n \notin I^N$, with $\zeta_N$ as in~\eqref{def:L_N}, and that $t \mapsto t/((t^2-1)^2+\mu)$ is decreasing for $t^2 \geq 1 + |\mu|$.
For the second column, for $v \in \ell^1_\nu$,
\[
  \|(\partial z)^* v\|_{X_1} = |(v,\partial_1 z)_2| + |(v,\partial_2 z)_2| \leq \|v\|_2 \bigl(\|\partial_1 z\|_1 + \|\partial_2 z\|_1\bigr) \leq 2\kappa\, r_2 \|v\|_{X_2},
\]
by Lemma~\ref{lem:derivative-estimate} and $\|v\|_2 \leq \|v\|_{X_2}$.
\end{proof}

\section{Results}
\label{sec:results}

\begin{figure}[H]
\centering
 \begin{minipage}{.5\linewidth}
  \centering\includegraphics[width=\linewidth]{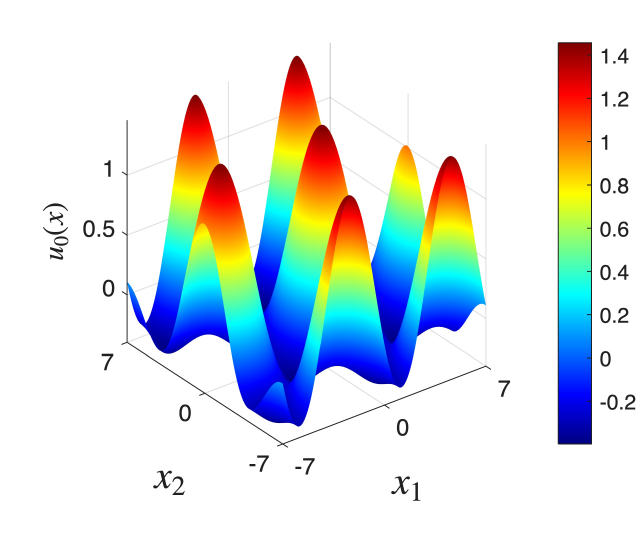}
 \end{minipage}%
 \begin{minipage}{.5\linewidth}
  \centering\includegraphics[width=\linewidth]{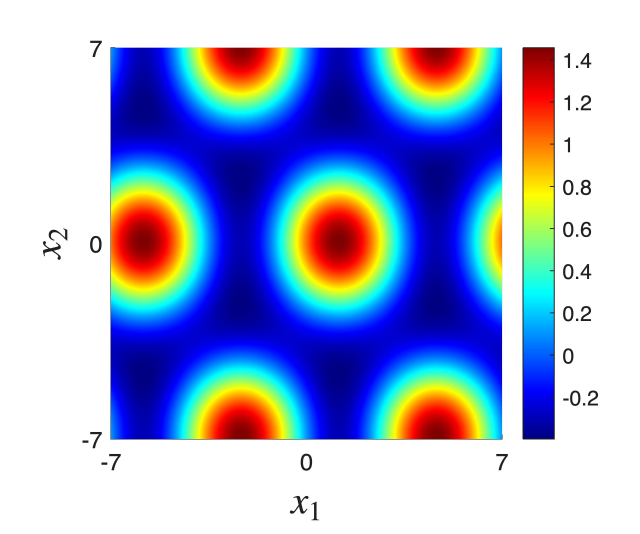}
 \end{minipage} \caption{Approximation of a stable solution in the 2D Swift Hohenberg equation with $\mu = 0.25, \gamma = 1.6$.}
 \label{fig:stable-solution}
 \end{figure}

\begin{figure}[H]
\centering
 \begin{minipage}{.5\linewidth}
  \centering\includegraphics[width=\linewidth]{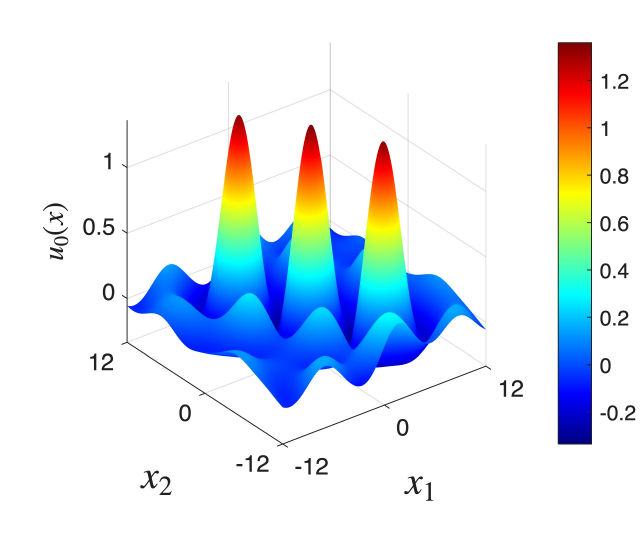}
 \end{minipage}%
 \begin{minipage}{.5\linewidth}
  \centering\includegraphics[width=\linewidth]{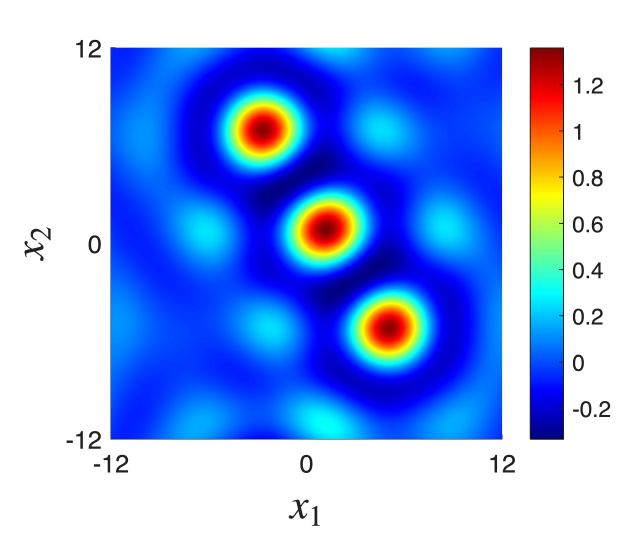}
 \end{minipage} \caption{Approximation of a stable solution in the 2D Swift Hohenberg equation with $\mu = 0.25, \gamma = 1.6$.}
 \label{fig:2stable-solution}
 \end{figure}

In this section, we present our computer-assisted proofs. Here stability is meant as in Section~\ref{sec:variational}: solutions of $\partial_t u = -F(u)$ on the same periodic domain $\om$ that start close to $\tilde{u}$ converge to a translate of $\tilde{u}$. The four examples show the range of outcomes: two stable patterns, and two unstable patterns with four and one unstable directions, respectively. In each case the number of unstable directions and the spectral gap come out of the same computation as the existence proof. In all four proofs below $N_i \geq \frac{d_i}{\pi}\sqrt{1 + |\mu|} - 1$, so the assumption $\zeta_N^2 \geq 1 + |\mu|$ of Lemma~\ref{lem:Zq} holds. 
\begin{theorem}[\bf The first stable solution]\label{th:stable one}
Let $\mu = 0.25, \gamma = 1.6$. Moreover, let $r \bydef (5 \times 10^{-8}, 5 \times 10^{-8})$. Then there exists a unique zero $\tilde{w}$ of $\overline{F}$ in $B_r(w_0) \subset X_{1} \times X_2$, and hence a stationary solution $\tilde{u}$ of the SH equation~\eqref{eq:swift-hohenberg-2d}, locally unique up to translations. The solution is stable. The approximation is shown in Figure \ref{fig:stable-solution}. Every nonzero eigenvalue $\lambda$ of $DF(\tilde{u})$ satisfies $|\lambda| \geq 0.0949$.
\end{theorem}
\begin{proof}
Choose $N_0 = (16,16), N = (10,10), d_1 = d_2 = 7, \nu_1 = \nu_2 = 1.1, \chi = \xi = (1,1)$. Then, we perform the full construction described in Section \ref{sec:numerical-aspects} to build $u_0$. We then define $w_0 \bydef (0,0,u_0)$. Next, we construct $\overline{A}$ using the approach described later in that section. Finally, using \cite{julia_blanco_symmetry_stab}, we choose $r = r^* \bydef (5 \times 10^{-8}, 5 \times 10^{-8})$ and define
{\small\begin{align}
    &H_1^{[1]}(r^*) \bydef 6.78 \times 10^{-6}, ~ H_1^{[2]} \bydef 6.98 \times 10^{-5}, ~ H_2^{[1]} \bydef 5.54 \times 10^{-4}, ~ H_2^{[2]} \bydef 2.97 \times 10^{-2} \\
    &Z_{Q,1}^{[1]}(r^*) \bydef 1.13 \times 10^{-7},~Z_{Q,1}^{[2]}(r^*) \bydef 1.06 \times 10^{-6},~ Z_{Q,2}^{[1]}(r^*) \bydef 3.02 \times 10^{-16},~Z_{Q,2}^{[2]}(r^*) \bydef 1.83 \times 10^{-6}
\end{align}}
and prove that these values satisfy Theorem \ref{th:radii-polynomial}. We also verify that $0 \notin \det \mathbf{M}$, with $\mathbf{M}$ as in~\eqref{eq:C-enclosure}, so that $M$ is invertible and Theorem~\ref{beta=0-2d} applies. Moreover, we have that
\begin{align}
    H_j^{[1]}(r^*) + Z_{Q,j}^{[1]}(r^*) + H_j^{[2]}(r^*) + Z^{[2]}_{Q,j}(r^*) < 1
\end{align}
so that \eqref{first-assumption-in-homotopy-lemma} and \eqref{second-assumption-in-homotopy-lemma} are satisfied. Then, using Gershgorin's theorem, we obtain that $\overline{A}^N$ has exactly two negative eigenvalues. The two negative discs and the positive disc closest to zero are approximately
\begin{align}
    D_1 \approx -3.114 \pm 8.30 \times 10^{-11}, D_2 \approx -2.633 \pm 8.30 \times 10^{-11}, D_3 \approx 0.00065 \pm 8.36 \times 10^{-11}.
\end{align}
According to Corollary \ref{lem:homotopy-arguments}, this shows that the solution is stable. 
\end{proof}
\begin{theorem}[\bf The second stable solution]
Let $\mu = 0.25, \gamma = 1.6$. Moreover, let $r \bydef (1 \times 10^{-10}, 1 \times 10^{-8})$. Then there exists a unique zero $\tilde{w}$ of $\overline{F}$ in $B_r(w_0) \subset X_{1} \times X_2$, and hence a stationary solution $\tilde{u}$ of the SH equation~\eqref{eq:swift-hohenberg-2d}, locally unique up to translations. The solution is stable. The approximation is shown in Figure \ref{fig:2stable-solution}. Every nonzero eigenvalue $\lambda$ of $DF(\tilde{u})$ satisfies $|\lambda| \geq 7.02 \times 10^{-3}$.
\end{theorem}
\begin{proof}
Choose $N_0 = (27,27), N = (20,20), d_1 = d_2 = 12, \nu_1 = \nu_2 = 1.09, \chi_1 = \chi_2 = 1$. The proof then follows in a similar way to that of Theorem \ref{th:stable one}. 
\end{proof}

\begin{figure}[H]
\centering
 \begin{minipage}{.5\linewidth}
  \centering\includegraphics[width=\linewidth]{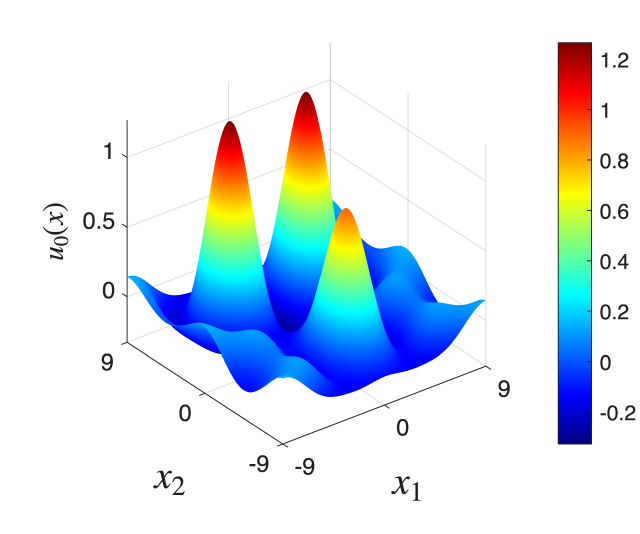}
 \end{minipage}%
 \begin{minipage}{.5\linewidth}
  \centering\includegraphics[width=\linewidth]{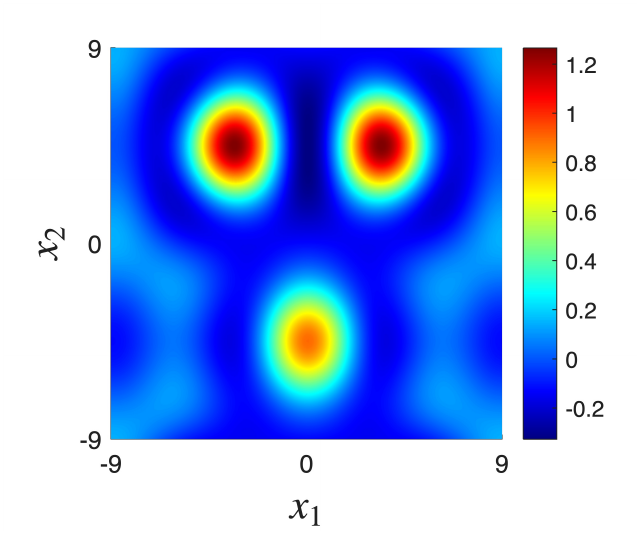}
 \end{minipage} \caption{Approximation of an unstable solution in the 2D Swift Hohenberg equation with $\mu = 0.25, \gamma = 1.6$.}
 \label{fig:unstable-solution}
 \end{figure}

\begin{figure}[H]
\centering
 \begin{minipage}{.5\linewidth}
  \centering\includegraphics[width=\linewidth]{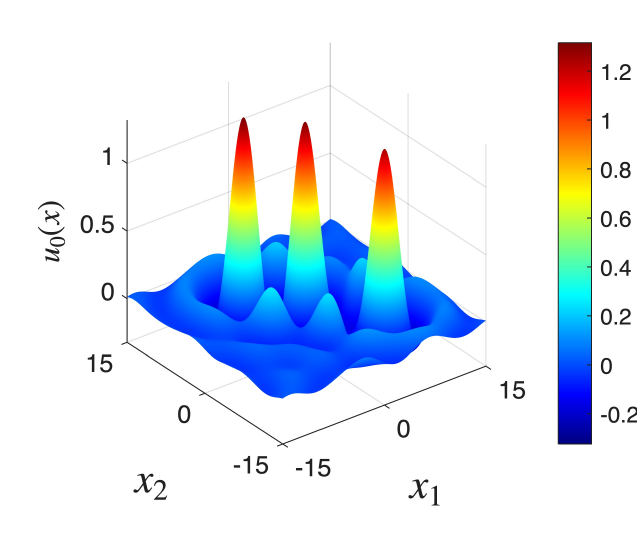}
 \end{minipage}%
 \begin{minipage}{.5\linewidth}
  \centering\includegraphics[width=\linewidth]{sh_sol_2unstable_2D_view2.pdf}
 \end{minipage} \caption{Approximation of an unstable solution in the 2D Swift Hohenberg equation with $\mu = 0.25, \gamma = 1.6$.}
 \label{fig:2unstable-solution}
 \end{figure}

 \begin{theorem}[\bf The first unstable solution]
Let $\mu = 0.25, \gamma = 1.6$. Moreover, let $r \bydef (1\times 10^{-9}, 6\times 10^{-8})$. Then there exists a unique zero $\tilde{w}$ of $\overline{F}$ in $B_r(w_0) \subset X_{1} \times X_2$, and hence a stationary solution $\tilde{u}$ of the SH equation~\eqref{eq:swift-hohenberg-2d}, locally unique up to translations. The solution is unstable with four unstable directions. The approximation is shown in Figure \ref{fig:unstable-solution}. Every nonzero eigenvalue $\lambda$ of $DF(\tilde{u})$ satisfies $|\lambda| \geq 1.22 \times 10^{-2}$.
\end{theorem}
\begin{proof}
Choose $N_0 = (18,18), N = (13,13), d_1 = d_2 = 9, \nu_1 = \nu_2 = 1.1, \chi_1 = \chi_2 = 1$. The proof then follows in a similar way to that of Theorem \ref{th:stable one}.
\end{proof}
  \begin{theorem}[\bf The second unstable solution]
Let $\mu = 0.25, \gamma = 1.6$. Moreover, let $r \bydef (4 \times 10^{-9}, 4 \times 10^{-9})$. Then there exists a unique zero $\tilde{w}$ of $\overline{F}$ in $B_r(w_0) \subset X_{1} \times X_2$, and hence a stationary solution $\tilde{u}$ of the SH equation~\eqref{eq:swift-hohenberg-2d}, locally unique up to translations. The solution is unstable with one unstable direction. The approximation is shown in Figure \ref{fig:2unstable-solution}. Every nonzero eigenvalue $\lambda$ of $DF(\tilde{u})$ satisfies $|\lambda| \geq 2.34 \times 10^{-3}$.
\end{theorem}
\begin{proof}
Choose $N_0 = (35,35), N = (27,27), d_1 = d_2 = 15, \nu_1 = \nu_2 = 1.08, \chi_1 = \chi_2 = 1$. The proof then follows in a similar way to that of Theorem \ref{th:stable one}.
\end{proof}

\section{Conclusion}
\label{sec:conclusion}

For variational problems with a continuous symmetry, we developed a method to prove the existence of a periodic solution and to determine its stability. Because of translation invariance, the solutions are not isolated. To obtain an isolated zero, we added phase conditions, with the unfolding parameters as Lagrange multipliers, which preserves the variational structure, see~\eqref{def:F_bar-2d}. Using the translation invariance of the energy, we then showed that the unfolding parameters vanish at a zero of the augmented system, provided the $2 \times 2$ matrix $M$ of inner products between the translation directions $\partial_j u_0$ of the approximate solution and $\partial_i \tilde{u}$ of the solution $\tilde{u}$ is invertible (Theorem~\ref{beta=0-2d}). Invertibility of $M$ means that no translation direction of $\tilde{u}$ is orthogonal to both translation directions of $u_0$, so that the phase conditions, which are built from $\partial_1 u_0$ and $\partial_2 u_0$, fix the position of $\tilde{u}$. Beyond the existence proof, this is the only condition needed. Although $\tilde{u}$ is not known explicitly, $M$ can be enclosed using the radius of the existence proof, so its invertibility can be checked with interval arithmetic. Determining the Morse index of $\tilde{u}$ and showing that the nullspace of the linearization $DF(\tilde{u})$ is spanned by $\partial_1 \tilde{u}$ and $\partial_2 \tilde{u}$ then requires no bounds beyond those of the existence proof (Remark~\ref{rem:border-u0}).

In addition, we obtained information on the eigenvalues of $DF(\tilde{u})$ beyond their count. For this we bordered $DF(\tilde{u})$ with its own translation directions $\partial_1 \tilde{u}$ and $\partial_2 \tilde{u}$, rather than with those of $u_0$ as in the augmented system. The approximate inverse of the existence proof is also an approximate inverse of the resulting operator $Q$ (Section~\ref{sec:stability-2D}), provided one additional bound $Z_Q$ is small enough, and its norm then gives a lower bound on the absolute values of the nonzero eigenvalues of $DF(\tilde{u})$ (Theorem~\ref{th:spectral-gap}). The bounds were computed for the Swift-Hohenberg equation in Section~\ref{sec:bounds-2D}, and the results are presented in Section~\ref{sec:results}.

A first limitation of our approach is the numerical cost. Since the solutions we prove do not have any symmetry, we cannot restrict to cosine series, so the matrix $\overline{A}^N$ has dimension $(2N_1+1)(2N_2+1)+2$, about four times that of a cosine-cosine expansion with the same $N_1, N_2$, and computing and storing it limits the size of the domain and the complexity of the patterns that can be treated. A second limitation is that the vanishing of the unfolding parameters relies on the variational structure. This is not a weakness of the argument: without variational structure, periodic patterns without reflection symmetry in general travel, and with $\partial_j u$ in place of $\partial_j u_0$ in the unfolding, the unfolding parameters become the velocity of the pattern.

We mention four extensions. First, the linear algebra in Sections~\ref{sec:beta-zero} and~\ref{sec:stability-2D} does not depend on the dimension, so the method works in any dimension, with one unfolding parameter for each translation direction.

Second, the method combines with the symmetry reduction of~\cite{JB_symmetries_1,JB_symmetries_2}. For some symmetry groups, the translations in certain directions commute with all the symmetries, so that the reduced problem still has a continuous symmetry. The present approach then applies with one unfolding parameter for each such direction. Stability refers to symmetric perturbations in that case.

Third, the argument is not specific to translations. For any other continuous symmetry of the energy, such as rotations or a phase symmetry, the derivative of the solution along the symmetry (for rotations $x_1\partial_2 \tilde{u} - x_2\partial_1 \tilde{u}$) lies in the nullspace and takes the role of $\partial_j \tilde{u}$ in the unfolding and in the bordered operators. Estimates such as Lemma~\ref{lem:derivative-estimate} must then be replaced by bounds for the derivative along the symmetry, in a function space suited to that symmetry. A natural application is to localized patterns in the plane. There the symmetries are the translations and the rotations, so three unfolding parameters are needed, and for $\mu > 0$ the essential spectrum of the linearization is bounded away from zero, so that the Morse index and the nullspace are again well defined. This would give proofs of localized patterns without the dihedral symmetry imposed in~\cite{symmetry_blanco_cadiot}.

Finally, the estimate behind Theorem~\ref{th:spectral-gap} contains more than a lower bound on the eigenvalues. The operators $Q(\tilde{w})^{-1}$ and $\overline{A}$ are self-adjoint and close in $\mathcal{B}(\ell^2)$, so their spectra are close~\cite[Theorem~V.4.10]{kato_1995}, and with the min-max principle their eigenvalues can be matched one by one. Since the spectrum of $\overline{A}$ is computable, this would give rigorous enclosures of the eigenvalues of $DF(\tilde{u})$ of smallest absolute value, such as the unstable eigenvalues of an unstable pattern and the slowest decay rate of a stable one. For this, the four eigenvalues $\pm\lambda_1,\pm\lambda_2$ of $Q(\tilde{w})$ caused by the border must first be identified. This can be done, for instance, with enclosures of $\lambda_1,\lambda_2$: their squares are the eigenvalues of the matrix $\Phi$ of Lemma~\ref{lem:bordered}, which can be enclosed in the same way as $M$, see~\eqref{eq:C-enclosure}.

\appendix
\section{Computing the Bounds}
\label{apen:computing-the-bounds}

We prove Lemmas~\ref{lem:Z2} and~\ref{lem:Z1N-computation}.

\subsection{The Bounds \texorpdfstring{$W_{22}^{[i]}$}{W22}}\label{apen:Z2}

We prove Lemma~\ref{lem:Z2}.
\begin{proof}
Let $0 < r \leq r^*$ and $w = w_0 + h \in B_r(w_0)$, and write $h = (\eta,v)$, so that $\|v\|_{X_2} \leq r_2 \leq r_2^*$. We show that the left-hand sides in Lemma~\ref{lem:Z2} are bounded by $W_{22}^{[i]} \|v\|_{X_2} \leq W_{22}^{[i]} r_2$.
Observe that
\begin{align}
    D\overline{F}(w_0 + h) - D\overline{F}(w_0) &= \begin{bmatrix}
         0 & 0 \\
        0 & DG(u_0 + v) - DG(u_0)
    \end{bmatrix} 
    \end{align}
Then, we have
    \begin{align}
        DG(u_0 + v) - DG(u_0) &= -2\gamma (\mathbb{u}_0 + \mathbb{v}) + 3(\mathbb{u}_0 + \mathbb{v})^2 + 2\gamma \mathbb{u}_0 - 3\mathbb{u}_0^2 = (-2\gamma I + 6\mathbb{u}_0)\mathbb{v} + 3\mathbb{v}^2.
\end{align}
We begin with $W_{22}^{[1]}$:
\begin{align}
    \|P_1 \overline{A}(D_u\overline{F}(w_0 + h) - D_u\overline{F}(w_0))\|_{\mathcal{B}(X_2,X_1)} &=\left\|P_1 \begin{bmatrix}
        0 & \overline{A}_{12}((-2\gamma I + 6\mathbb{u}_0)\mathbb{v} + 3\mathbb{v}^2) \\
        0 & \overline{A}_{22} ((-2\gamma I + 6\mathbb{u}_0)\mathbb{v} + 3\mathbb{v}^2)
    \end{bmatrix}\right\|_{\mathcal{B}(X_2,X_1)} \\
    &=\left\|\overline{A}_{12}((-2\gamma I + 6\mathbb{u}_0)\mathbb{v} + 3\mathbb{v}^2)\right\|_{\mathcal{B}(X_2,X_1)} \\
    &\hspace{-2cm}\leq \left\| \overline{A}_{12}(-2\gamma I + 6\mathbb{u}_0)\mathbb{v}\right\|_{\mathcal{B}(X_2,X_1)} + 3\left\| \overline{A}_{12}\mathbb{v}^2
    \right\|_{\mathcal{B}(X_2,X_1)} \\
    &\hspace{-2.5cm}\leq \left\| \overline{A}_{12}\mathbb{q}\right\|_{\mathcal{B}(X_2,X_1)} \|\mathbb{v}\|_{\mathcal{B}(X_2)} + 3\| \overline{A}_{12}\|_{\mathcal{B}(X_2,X_1)} \|\mathbb{v}^2\|_{\mathcal{B}(X_2)}\label{before-banach-algebra}
\end{align}
where $q = (q_n)_{n \in \mathbb{Z}^2} = (-2\gamma \delta_n + 6a_n)_{n \in \mathbb{Z}^2}$, and $\delta_n$ is the Kronecker delta. Since $\overline{A}_{12} = \overline{A}_{12}\pi^{\le N}$ and $q$ has its Fourier coefficients in $I^{N_0}$, we have $\overline{A}_{12}\mathbb{q} = \overline{A}_{12}\mathbb{q}\pi^{\le N+N_0}$, and similarly for $\overline{A}_{22}^N$.
Using that $\ell^1_\nu$ is a Banach algebra, observe that
\begin{align}
    &\|\mathbb{v}\|_{\mathcal{B}(\ell^1_\nu)} = \sup_{z \in \ell^1_\nu} \frac{\|\mathbb{v}z\|_{\ell^1_\nu}}{\|z\|_{\ell^1_\nu}} = \sup_{z \in \ell^1_\nu} \frac{\|v*z\|_{\ell^1_\nu}}{\|z\|_{\ell^1_\nu}} \leq\sup_{z \in \ell^1_\nu} \frac{\|v\|_{\ell^1_\nu}\|z\|_{\ell^1_\nu}}{\|z\|_{\ell^1_\nu}} = \|v\|_{\ell^1_\nu} \\
    &\|\mathbb{v}^2\|_{\mathcal{B}(\ell^1_\nu)} = \sup_{z \in \ell^1_\nu} \frac{\|\mathbb{v}^2z\|_{\ell^1_\nu}}{\|z\|_{\ell^1_\nu}} = \sup_{z \in \ell^1_\nu} \frac{\|v*v*z\|_{\ell^1_\nu}}{\|z\|_{\ell^1_\nu}} \leq\sup_{z \in \ell^1_\nu} \frac{\|v\|_{\ell^1_\nu}^2\|z\|_{\ell^1_\nu}}{\|z\|_{\ell^1_\nu}} = \|v\|_{\ell^1_\nu}^2 \leq r_2^* \|v\|_{\ell^1_\nu}.
\end{align}
Returning to \eqref{before-banach-algebra}, we have
\begin{align}
    \|P_1\overline{A}(D_u\overline{F}(w_0 + h) - D_u\overline{F}(w_0))\|_{\mathcal{B}(X_2,X_1)} &\leq \left\| \overline{A}_{12}\mathbb{q}
    \right\|_{\mathcal{B}(X_2,X_1)} \|v\|_{X_2} + 3\|\overline{A}_{12}\|_{\mathcal{B}(X_2,X_1)} r_2^* \|v\|_{X_2} = W_{22}^{[1]} \|v\|_{X_2}
    \end{align}
as desired. 
For $W_{22}^{[2]}$, the same steps give
\begin{align}
    \|P_2 \overline{A}(D_u \overline{F}(w_0 + h) - D_u \overline{F}(w_0))\|_{\mathcal{B}(X_2)} \leq \bigl(\|\overline{A}_{22} \mathbb{q}\|_{\mathcal{B}(X_2)} + 3\|\overline{A}_{22}\|_{\mathcal{B}(X_2)} r_2^*\bigr) \|v\|_{X_2}.\label{W_2^2-term-before-handling-A_22}
\end{align}
This term is handled differently since $\|\pi^{>N} \overline{A}_{22}\|_{\mathcal{B}(X_2)} \leq \frac{1}{L_N}$. We write
\begin{align}
    \|\overline{A}_{22} \mathbb{q}\|_{\mathcal{B}(X_2)} &\leq \left\| \overline{A}_{22}^N\mathbb{q}\right\|_{\mathcal{B}(X_2)} + \left\| \pi^{>N}\overline{A}_{22}\mathbb{q}\right\|_{\mathcal{B}(X_2)}\\&\leq \left\|\overline{A}_{22}^N\mathbb{q}\pi^{\le N+N_0}\right\|_{\mathcal{B}(X_2)}+\|\pi^{>N} \overline{A}_{22}\|_{\mathcal{B}(X_2)} \| \mathbb{q}\|_{\mathcal{B}(X_2)} \\
    &\leq \left\|\overline{A}_{22}^N\mathbb{q}\pi^{\le N+N_0}
    \right\|_{\mathcal{B}(X_2)}+\frac{1}{L_N} \|q\|_{\ell^1_\nu}
\end{align}
where we used the fact that 
\begin{align}
    \|\mathbb{q}\|_{\mathcal{B}(\ell^1_\nu)} &= \sup_{z \in \ell^1_\nu} \frac{\|\mathbb{q}z\|_{\ell^1_\nu}}{\|z\|_{\ell^1_\nu}} =\sup_{z \in \ell^1_\nu} \frac{\|qz\|_{\ell^1_\nu}}{\|z\|_{\ell^1_\nu}} \leq \sup_{z \in \ell^1_\nu} \frac{\|q\|_{\ell^1_\nu}\|z\|_{\ell^1_\nu}}{\|z\|_{\ell^1_\nu}} = \|q\|_{\ell^1_\nu}.\label{first-estimate-W_2^2}
\end{align}
Returning to \eqref{W_2^2-term-before-handling-A_22}, we examine the second term and estimate
\begin{align}
    \|\overline{A}_{22}\|_{\mathcal{B}(X_2)} \leq \|\overline{A}^N_{22}\|_{\mathcal{B}(X_2)} + \|\pi^{>N} \overline{A}_{22}\|_{\mathcal{B}(X_2)} \leq \|\overline{A}^N_{22}\|_{\mathcal{B}(X_2)} + \frac{1}{L_N}.\label{second-estimate-W_2^2}
\end{align}
Combining \eqref{first-estimate-W_2^2} and \eqref{second-estimate-W_2^2}, we obtain
\begin{align}
    &\|P_2\overline{A}(D_u\overline{F}(w_0 + h) - D_u\overline{F}(w_0))\|_{\mathcal{B}(X_2)} \\
    &\leq \left(\left\|\overline{A}_{22}^N\mathbb{q}\pi^{\le N+N_0}
    \right\|_{\mathcal{B}(X_2)}+\frac{1}{L_N} \|q\|_{\ell^1_\nu} + 3\left(\|\overline{A}_{22}^N\|_{\mathcal{B}(X_2)} + \frac{1}{L_N}\right)r_2^*\right) \|v\|_{X_2} \\
    &= W_{22}^{[2]} \|v\|_{X_2}
\end{align}
as desired.    
\end{proof}

\subsection{The Bounds \texorpdfstring{$Z_j^{[i],N}$}{Zj\^iN}}\label{apen:Z1N}
We prove Lemma~\ref{lem:Z1N-computation}. There are multiple approaches for computing each of these bounds. We choose an approach that relies on the computer. Indeed, by using non-square operators, we are able to allow the computer to do more of the work for us and perform less analysis by hand. While this approach does require additional memory, it was sufficient for our purposes.
\begin{proof}
Note that $Z_{1}^{[1],N}$ can be computed directly as it is an operator in $X_1$. For $Z_2^{[1],N}$, we split off the modes outside $I^{2N}$:
\begin{align}
    &\|\overline{A}_{11}(\partial u_0^N)^* + \overline{A}_{12} (L + \mathbb{v}_0^N)\|_{\mathcal{B}(X_2,X_1)} \\
    &\leq \|(\overline{A}_{11}(\partial u_0^N)^* + \overline{A}_{12} (L + \mathbb{v}_0^N))\pi^{\le 2N}\|_{\mathcal{B}(X_2,X_1)} + \|(\overline{A}_{11}(\partial u_0^N)^* + \overline{A}_{12} (L + \mathbb{v}_0^N))\pi^{>2N}\|_{\mathcal{B}(X_2,X_1)} .
\end{align}
The second term vanishes, since $\overline{A}_{12} = \overline{A}_{12}\pi^{\le N}$, $(\partial u_0^N)^* = (\partial u_0^N)^*\pi^{\le N}$, and
\begin{align}
    \pi^{\le N}(L + \mathbb{v}_0^N)\pi^{>2N} = 0,\label{getting-0-in-Z1N-proof}
\end{align}
because $L$ is diagonal and multiplication by $v_0^N$ shifts Fourier modes by at most $N$. Hence, we define
\begin{align}
    Z_{2}^{[1],N} \bydef \|(\overline{A}_{11}(\partial u_0^N)^* + \overline{A}_{12} (L + \mathbb{v}_0^N))\pi^{\le 2N}\|_{\mathcal{B}(X_2,X_1)}.
\end{align}
For $Z_1^{[2],N}$ we again start from the definition.
\begin{align}
    \|\overline{A}_{22} \partial u_0^N\|_{\mathcal{B}(X_1,X_2)} &= \|\overline{A}_{22}^N \partial u_0^N\|_{\mathcal{B}(X_1,X_2)} + \|\pi^{>N} \overline{A}_{22} \partial u_0^N\|_{\mathcal{B}(X_1,X_2)} \\
    &\leq \|\overline{A}_{22}^N \partial u_0^N\|_{\mathcal{B}(X_1,X_2)} + \|\pi^{>N} \overline{A}_{22}\|_{\mathcal{B}(X_2)} \|\pi^{>N} \partial u_0^N\|_{\mathcal{B}(X_1,X_2)} \\
    &= \|\overline{A}_{22}^N \partial u_0^N\|_{\mathcal{B}(X_1,X_2)} \bydef Z_1^{[2],N}
\end{align}
where the last term was $0$ due to $\pi^{>N} \partial u_0^N = 0$. Let us finally consider $Z_2^{[2],N}$. Then, observe that
\begin{align}
    &I - \overline{A}_{21} (\partial u_0^N)^* - \overline{A}_{22} (L + \mathbb{v}_0^N) \\
    &= \pi^{\le N}(I - \overline{A}_{21} (\partial u_0^N)^* - \overline{A}_{22} (L + \mathbb{v}_0^N))\pi^{\le 2N} + \pi^{\le N}(I - \overline{A}_{21} (\partial u_0^N)^* - \overline{A}_{22} (L + \mathbb{v}_0^N))\pi^{>2N}\\
    &\hspace{+1cm} + \pi^{>N}(I - \overline{A}_{21} (\partial u_0^N)^* - \overline{A}_{22} (L + \mathbb{v}_0^N)) \\
    &= \pi^{\le N}( I - \overline{A}_{21} (\partial u_0^N)^* - \overline{A}_{22} (L + \mathbb{v}_0^N)) \pi^{\le 2N} - \overline{A}_{22}^N (L + \mathbb{v}_0^N) \pi^{>2N} + \pi^{> N} - \pi^{>N} \overline{A}_{22} (L + \mathbb{v}_0^N).
\end{align}
The first term is the norm of a finite matrix, which we enclose with interval arithmetic. We define
\begin{align}
    Z_{2}^{[2],N,1} \bydef \|\pi^{\le N}( I - \overline{A}_{21} (\partial u_0^N)^* - \overline{A}_{22} (L + \mathbb{v}_0^N)) \pi^{\le 2N}\|_{\mathcal{B}(X_2)}.
\end{align}
For the second term, notice that
\begin{align}
    \|\overline{A}_{22}^N (L + \mathbb{v}_0^N) \pi^{>2N}\|_{\mathcal{B}(X_2)} &\leq \|\overline{A}_{22}^N\|_{\mathcal{B}(X_2)} \|\pi^{\le N} (L + \mathbb{v}_0^N) \pi^{>2N}\|_{\mathcal{B}(X_2)} = 0
\end{align}
where we followed the same steps as in \eqref{getting-0-in-Z1N-proof}. Hence, this term is zero. For the third term,
{\footnotesize\begin{align}
    \|\pi^{> N} - \pi^{>N} \overline{A}_{22} (L + \mathbb{v}_0^N)\|_{\mathcal{B}(X_2)} =\|\pi^{> N} - \pi^{>N} \overline{A}_{22}L - \pi^{> N} \overline{A}_{22} \mathbb{v}_0^N\|_{\mathcal{B}(X_2)}
    &= \|\pi^{> N} - \pi^{>N} - \pi^{> N} \overline{A}_{22} \mathbb{v}_0^N\|_{\mathcal{B}(X_2)} 
    \end{align}}
where we used the fact that $\pi^{>N} \overline{A}_{22} = \pi^{>N} L^{-1}$. Continuing, 
    \begin{align}
    \|\pi^{> N} - \pi^{>N} - \pi^{> N} \overline{A}_{22} \mathbb{v}_0^N\|_{\mathcal{B}(X_2)}  &= \|\pi^{>N} \overline{A}_{22}\mathbb{v}_0^N \|_{\mathcal{B}(X_2)} \\
    &\leq \|\pi^{>N} \overline{A}_{22}\|_{\mathcal{B}(X_2)} \|\pi^{>N} \mathbb{v}_0^N \|_{\mathcal{B}(X_2)} \\
    &\leq \frac{1}{L_N} \|\mathbb{v}_0^N\|_{\mathcal{B}(X_2)} \\
    &\leq \frac{1}{L_N} \|v_0^N\|_{\ell^1_\nu} \bydef Z_{2}^{[2],N,2}
\end{align}
where we used similar steps to those used to compute the $Z_{2}^{[i],\infty}$ bounds for $i = 1,2$. 
Hence, we define
\begin{align}
    Z_{2}^{[2],N} \bydef Z_{2}^{[2],N,1} + Z_{2}^{[2],N,2}
\end{align}
and obtain the final result.
\end{proof}

\section{Counting negative eigenvalues}
\label{apen:homotopy}

We prove the equality $\mathcal{N}(Q(\tilde{w})) = \mathcal{N}(\overline{A}^N)$ of Proposition~\ref{prop:homotopy}. We write $Q \bydef Q(\tilde{w})$, and we recall from Theorem~\ref{th:radii-polynomial-q} that $\|I - \overline{A}Q\|_\xi \leq \delta < 1$, with $\delta$ as defined in~\eqref{first-assumption-in-homotopy-lemma}. In particular $\overline{A}$ is surjective onto $X_1 \times X_2$, so the matrix $\overline{A}^N$ is invertible, and $\overline{A}^{-1} = (\overline{A}^N)^{-1} \oplus \pi^{>N} L$.

The homotopy $(1-s)\overline{A} + sQ^{-1}$, the analogue of the one used in~\cite{JB_symmetries_2}, passes through compact operators. Their spectrum accumulates at $0$, so a negative eigenvalue might tend to $0$ and disappear into this accumulation point, without $0$ ever being an eigenvalue. We avoid this by interpolating between unbounded operators with the same leading part $L$, whose eigenvalues accumulate only at $+\infty$. Let
\[
  E \bydef Q - \overline{A}^{-1}, \qquad Q_s \bydef \overline{A}^{-1} + sE, \quad s \in [0,1],
\]
so that $Q_0 = \overline{A}^{-1}$ and $Q_1 = Q$. We prove that $\mathcal{N}(Q_s)$ does not depend on $s$. We first establish that $E$ is bounded. It follows that the eigenvalues of $Q_s$ move continuously in $s$, so that $\mathcal{N}(Q_s)$ can only change at a value of $s$ for which $Q_s$ is not injective. Finally we show that $Q_s$ is injective for every~$s$, using $\|I - \overline{A}Q\|_\xi \leq \delta < 1$.

\emph{$E$ is bounded.} We also write $L$ for the operator $(\beta,\phi) \mapsto (0,L\phi)$ on $\mathbb{C}^2 \times \ell^2$. Then $Q - L$ consists of the convolution $V$ of Lemma~\ref{lem:regularity}, which is bounded on $\ell^2$, and of the border terms built from $\partial_j \tilde{u} \in \ell^2$, so it is bounded. Moreover, $\overline{A}^{-1} - L$ acts only on the finitely many coordinates $\beta$ and $\phi_n$ with $n \in I^N$, so it is a matrix. Hence $E = (Q - L) - (\overline{A}^{-1} - L)$ is a bounded self-adjoint operator on $\mathbb{C}^2 \times \ell^2$.

\emph{Eigenvalues move continuously.} Since $E$ is bounded and self-adjoint, each $Q_s$ is, like $\overline{A}^{-1}$, self-adjoint and bounded below as a quadratic form, with compact resolvent, and all $Q_s$ have the same domain $\mathbb{C}^2 \times \mathrm{Dom}(L)$, where $\mathrm{Dom}(L)$ consists of the $u \in \ell^2$ with $(L_{nn} u_n)_n \in \ell^2$. Hence the spectrum of $Q_s$ consists of eigenvalues $\mu_1(s) \leq \mu_2(s) \leq \dots$, counted with multiplicity, with $\mu_k(s) \to \infty$ as $k \to \infty$. By the min-max principle, see e.g.~\cite[Theorem~XIII.1]{reed_simon_1978},
\[
  \mu_k(s) = \min_{\dim S = k} \, \max_{x \in S, \, \|x\| = 1} (Q_s x,x),
\]
where $S$ runs over the $k$-dimensional subspaces of the domain. For $\|x\| = 1$ we have $(Q_s x,x) = (Q_t x,x) + (s-t)(Ex,x)$ with $|(Ex,x)| \leq \|E\|_{\mathcal{B}(\ell^2)}$, so the quantity inside the min-max changes by at most $|s-t|\,\|E\|_{\mathcal{B}(\ell^2)}$ when $t$ is replaced by $s$, and therefore
\[
  |\mu_k(s) - \mu_k(t)| \leq |s-t|\,\|E\|_{\mathcal{B}(\ell^2)} .
\]
Hence each $\mu_k$ is continuous in $s$. If no $\mu_k$ vanishes on $[0,1]$, then each $\mu_k$ keeps its sign, so $\mathcal{N}(Q_s)$ can only change at a value of $s$ for which some $\mu_k(s) = 0$. Since the spectrum of $Q_s$ consists of eigenvalues, $Q_s$ is then not injective. For more general results on the continuity of eigenvalues under perturbations, see~\cite[Section~IV.3.5]{kato_1995}.

\emph{$Q_s$ is injective.} We have $\overline{A}Q_s = I + s(\overline{A}Q - I)$, so
\[
  \|I - \overline{A}Q_s\|_\xi = s\,\|I - \overline{A}Q\|_\xi \leq s\delta < 1 ,
\]
and $Q_s$ is injective on $X_1 \times X_2$. If $z = (\beta,\phi) \in \mathbb{C}^2 \times \ell^2$ satisfies $Q_s z = 0$, then, since $\overline{A}^{-1}$ acts as $L$ on the modes outside $I^N$, we have $\pi^{>N}\bigl(L\phi + s(V\phi + \beta_1\partial_1\tilde{u} + \beta_2\partial_2\tilde{u})\bigr) = 0$. Arguing as in the proof of Lemma~\ref{lem:regularity}, but with $sV$ instead of $V$ and with the additional term $s(\beta_1\partial_1\tilde{u} + \beta_2\partial_2\tilde{u})$, which is harmless since $L^{-1}\partial_j\tilde{u} \in \ell^1_\nu$, we conclude that $\phi \in \ell^1_\nu$. Hence $z \in X_1 \times X_2$ and $z = 0$.

We conclude that $\mathcal{N}(Q) = \mathcal{N}(Q_1) = \mathcal{N}(Q_0) = \mathcal{N}(\overline{A}^{-1})$. Finally $\mathcal{N}(\overline{A}^{-1}) = \mathcal{N}(\overline{A}^N)$, since $\pi^{>N} L$ is positive and the eigenvalues of $(\overline{A}^N)^{-1}$ are the reciprocals of those of $\overline{A}^N$.

\begin{remark}\label{rem:homotopy-DFbar}
The argument above, with $Q$ replaced by $D\overline{F}(\tilde{w})$, also gives $\mathcal{N}(D\overline{F}(\tilde{w})) = \mathcal{N}(\overline{A}^N)$, which is used in Remark~\ref{rem:border-u0}. The operator $D\overline{F}(\tilde{w})$ is self-adjoint and differs from $Q$ only in the border, so $D\overline{F}(\tilde{w}) - \overline{A}^{-1}$ is bounded. The estimate~\eqref{second-assumption-in-homotopy-lemma} takes the place of~\eqref{first-assumption-in-homotopy-lemma}, with $\|\cdot\|_\chi$ in place of $\|\cdot\|_\xi$. In the regularity argument the border term becomes $s(\beta_1\partial_1 u_0 + \beta_2\partial_2 u_0)$, which has only finitely many nonzero Fourier coefficients.
\end{remark}

\newpage
\nocite{julia_interval}
\bibliographystyle{abbrv}
\bibliography{biblio}
\end{document}